\documentclass[]{article}
\usepackage[latin9]{inputenc}
\usepackage{float}
\usepackage{amsthm}
\usepackage{amsmath}
\usepackage{amssymb}
\usepackage{graphicx}

\usepackage[title]{appendix}

\makeatletter

\theoremstyle{plain}
\newtheorem{thm}{\protect\theoremname}
  \theoremstyle{plain}
  \newtheorem{prop}{\protect\propositionname}
  \theoremstyle{definition}
  \newtheorem{defn}{\protect\definitionname}
 \theoremstyle{definition}
  \newtheorem{example}{\protect\examplename}

\makeatother

  \providecommand{\definitionname}{Definition}
  \providecommand{\examplename}{Example}
  \providecommand{\propositionname}{Proposition}
\providecommand{\theoremname}{Theorem}

\begin{document}

\title{Shearless transport barriers in unsteady\\
two-dimensional flows and maps\footnote{Submitted to Physica D}}

\author{Mohammad Farazmand$^{1,2}$, Daniel Blazevski$^{2}$, George Haller$^{2}$%
\thanks{Corresponding author email: georgehaller@ethz.ch}\\
\ \\
{\small
$^{1}$Department of Mathematics, ETH Zurich, Rämistrasse 101, 8092 Zurich Switzerland}
\\
{\small $^{2}$Institute of Mechanical Systems, Department of Mechanical and
Process Engineering}\\
{\small ETH Zurich, Tannenstrasse 3, 8092 Zurich, Switzerland}
}

\maketitle

\begin{abstract}
We develop a variational principle that extends the notion of a shearless
transport barrier from steady to general unsteady two-dimensional
flows and maps defined over a finite time interval. This principle
reveals that hyperbolic Lagrangian Coherent Structures (LCSs) and
parabolic LCSs (or jet cores) are the two main types of shearless
barriers in unsteady flows. Based on the boundary conditions they
satisfy, parabolic barriers are found to be more observable and robust
than hyperbolic barriers, confirming widespread numerical observations.
Both types of barriers are special null-geodesics of an appropriate
Lorentzian metric derived from the Cauchy--Green strain tensor. Using
this fact, we devise an algorithm for the automated computation of
parabolic barriers. We illustrate our detection method on steady and
unsteady non-twist maps and on the aperiodically forced Bickley jet.
\end{abstract}



\section{Introduction}

A shearless transport barrier in two dimensions is generally defined
as a member of a closed invariant curve family whose frequency admits
a local extremum within the family. This definition ties shearless
barriers fundamentally to recurrent (i.e., steady, periodic or quasiperiodic)
flows where the necessary frequencies are well-defined. Here we extend
the notion of a shearless transport barrier to two-dimensional flows
and maps with general time-dependence.

\begin{figure}[H]
\begin{center}
\includegraphics[width=0.6\textwidth]{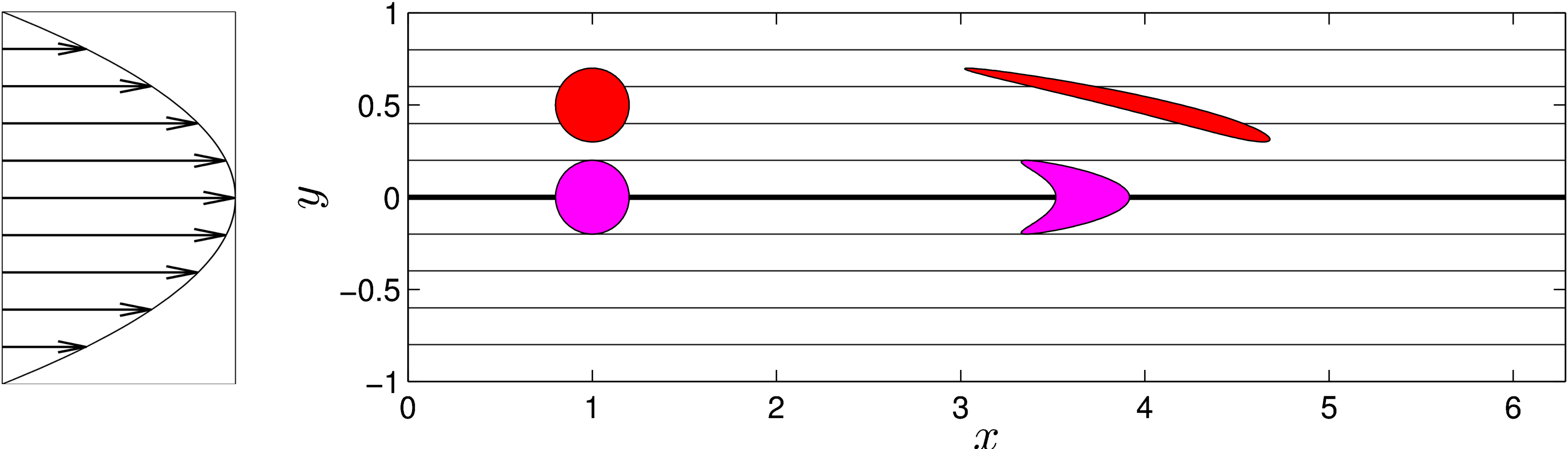}
\end{center}
\caption{Left: The velocity profile of the steady flow (\ref{eqn:steady_par_intro}) for $u(y)=1-y^{2}$. 
Right: Streamlines for the same flow.
The thick line at $y=0$ marks the shearless streamline that acts
as a jet core. The tracer disk located on the shearless line (magenta
circle) deforms into a blunt arrow shape symmetrically under advection
to time $t=9$. The tracer disk located away from the shearless line
(red circle) has a markedly different deformation pattern.}
\label{fig:canShearJet_tracers} 
\end{figure}

In steady and time-periodic problems of fluid dynamics and plasma
physics, shearless (or non-twist) barriers have been found to be particularly
robust inhibitors of phase space transport \cite{Diego_NT,Javier_PRL,Morrison_PRE,samelson_jet}.
For illustration, consider a steady, parallel shear flow
\begin{align}
\dot{x} & =u(y),\qquad u^{\prime}(y_{0})=0.\label{eqn:steady_par_intro}\\
\dot{y} & =0,\nonumber 
\end{align}
 on a domain periodic in $x$. The $y=y_{0}$ line marks a jet core,
whose impact on tracer patterns is shown in Fig. \ref{fig:canShearJet_tracers}
in a particular example with $y_{0}=0$. Note the unique material
signature of the shearless barrier, deforming the tracer blob initialized
along it into a boomerang-shaped pattern, By contrast, another tracer
blob simply stretches under shear.

The flow \eqref{eqn:steady_par_intro} is an idealized model of the
velocity field inside atmospheric or oceanic zonal jets, or helical
magnetic field lines in a tokamak \cite{Abdullaev}. As a dynamical
system, \eqref{eqn:steady_par_intro} represents an integrable system
with the Hamiltonian $H(y)=\int_{0}^{y}u(\eta)d\eta$. Its horizontal
trajectories along which the Eulerian shear $u'(y)$ vanishes are
referred to as shearless barriers. Along these barriers, $H''(y_{0})=0$
holds, thus the circle $y=y_{0}$ does not satisfy the twist condition
of classic KAM theory \cite{Arnold}.

Yet numerical studies of \cite{Diego_NT,Javier_PRL,Morrison_PRE,wurm05}
show that such barriers are more robust under steady or time-periodic
perturbations than any other nearby KAM tori. Related theoretical
results for two-dimensional maps were given in \cite{Rafael_NT_siam}.
More recently, degenerate tori for steady $3$D maps were considered
in \cite{Mezic_3D_deg}. In addition, a general \emph{a posteriori}
result on non-twist tori of arbitrary dimension that are potentially
far from integrable has been obtained by \cite{Rafael_NT_mem}. However,
no general theory of shearless transport barriers for unsteady flows
has been established.

The need for such a general theory of unsteady shearless barriers
clearly exists. In plasma physics, computational and experimental studies
suggest that shearless barriers 
enhance the confinement of plasma in magnetic fusion devices \cite{JET, Lorenzini_Nat_Phys, Lorenzini_Nuc_Fus, Rev_shear},
which generate turbulent velocity fields with general time dependence.  In this context, a description of
shearless barriers is either understood in models for steady magnetic fields \cite{Rev_shear} or inferred from scalar quantities (e.g. temperature, density) in more complex unsteady scenarios \cite{JET, Lorenzini_Nat_Phys, Lorenzini_Nuc_Fus}.

In fluid dynamics, shearless barriers are of interest in the context
of zonal jets. Rossby waves are the best known and most robust transport
barriers in geophysical flows \cite{Dfluid,Javier_bickley,Swinney},
yet only recent work attempts to their attendant unsteady jet cores
in the Lagrangian fame of an unsteady flow. The method put forward
in \cite{Javier_shear} seeks such Lagrangian shearless barriers as trenches
of the finite-time Lyapunov exponent (FTLE) field. However, just as
the examples in \cite{var_theory} show that FTLE ridges do not necessarily
correspond to hyperbolic Lagrangian structures, FTLE trenches may
also fail to mark zonal jet cores (see Example \ref{ex:FTLE-Trench}
in Section \ref{sub:Parabolic-barriers} below).

Here we develop a variational principle for shearless barriers as
centerpieces of material strips showing no leading order variation
in Lagrangian shear. This variational principle shows that shearless
barriers are composed of tensorlines of the right Cauchy--Green strain
tensor associated with the flow map. Most stretching or contracting
Cauchy--Green tensorlines have previously been identified as best
candidates for hyperbolic Lagrangian Coherent Structures (LCSs) \cite{geo_theory,stretchlines},
but no underlying global variational principle has been known to which
they would be solutions. The present work, therefore, also advances
the theory of hyperbolic LCS, establishing them as shearless transport
barriers under fixed (Dirichlet-type) boundary conditions. 

Our main result is that parabolic transport barriers (jet cores) are
also solutions of the same shearless Lagrangian variational principle,
satisfying variable-endpoint boundary conditions. They are formed
by minimally hyperbolic, structurally stable chains of tensorlines
that connect singularities of the Cauchy--Green strain tensor field.
We develop and test a numerical procedure that detects such tensorline
chains, thereby finding generalized Lagrangian jet cores in an arbitrary,
two-dimensional unsteady flow field in an automated fashion.

\section{Notation and definitions\label{sec:Notation}}

Let $v(x,t)$ denote a two-dimensional velocity field, with $x$ labeling
positions in a two-dimensional region $U$, and with $t$ referring
to time. Fluid trajectories generated by this velocity field satisfy
the differential equation 
\begin{equation}
\dot{x}=v(x,t),\label{eq:uxt}
\end{equation}
whose solutions are denoted by $x(t;t_{0},x_{0}$), with $x_{0}$
referring to the initial position at time $t_{0}$. The evolution
of fluid elements is described by the flow map 
\begin{equation}
F_{t_{0}}^{t}(x_{0}):=x(t;t_{0},x_{0}),
\end{equation}
which takes any initial position $x_{0}$ to its current position
at time $t$. 

Lagrangian strain in the flow is often characterized by the right
Cauchy--Green strain tensor field $C(x_{0})=\left[\nabla F_{t_{0}}^{t}(x_{0})\right]^{T}\nabla F_{t_{0}}^{t}(x_{0})$,
whose eigenvalues $\lambda_{i}(x_{0})$ and eigenvectors $\xi_{i}(x_{0})$
satisfy
\[
C\xi_{i}=\lambda_{i}\xi_{i},\quad\left|\xi_{i}\right|=1,\quad i=1,2;\qquad0<\lambda_{1}\leq\lambda_{2},\qquad\xi_{1}\perp\xi_{2}.
\]
 The tensor $C$, as well as its eigenvalues and eigenvectors, depend
on the choice of the times $t$ and $t_{0}$, but we suppress this
dependence for notational simplicity.

\section{Stability of material lines\label{sec:stability}}

\begin{figure}[H]
\begin{center}
\includegraphics[width=0.6\textwidth]{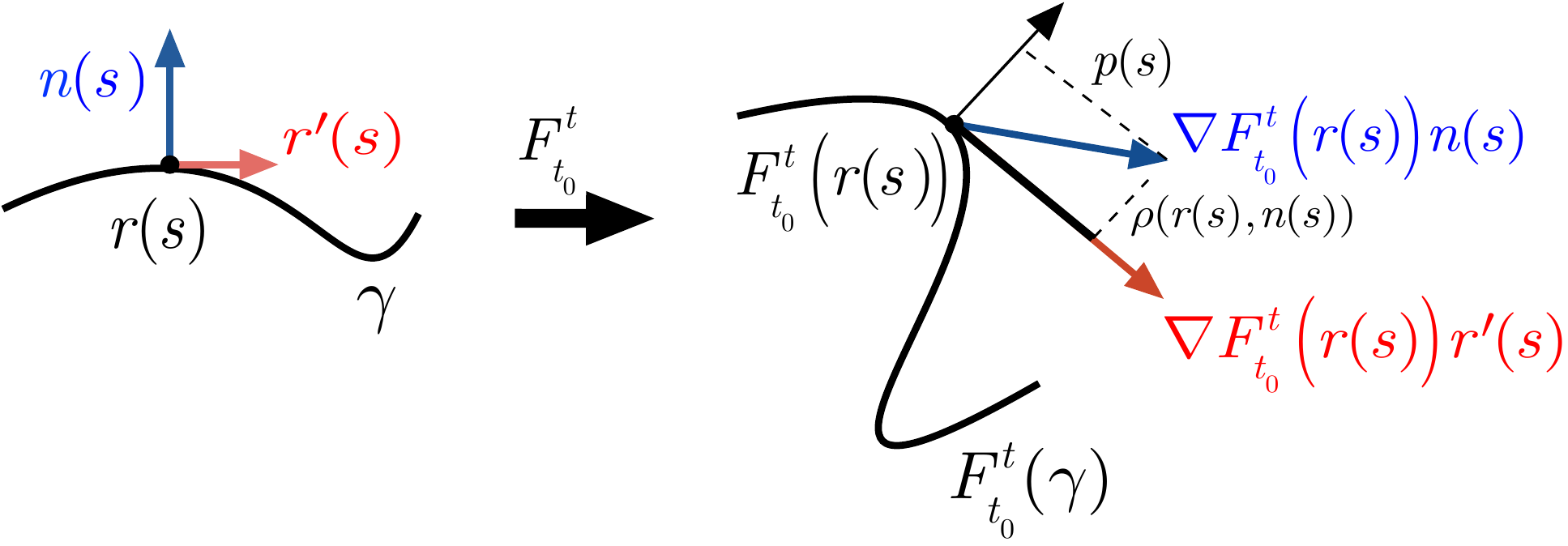}
\end{center}
\caption{The evolution of a unit normal vector $n(s)$ of a material line $\gamma$
under the linearized flow map $\nabla F_{t_{0}}^{t}$.}
\label{fig:rho}
\end{figure}

Consider a material line (i.e., a smooth curve of initial conditions)
$\gamma$ at time $t_{0}$, parametrized as $r(s)$ with $s\in[0,\sigma]$.
If $n(s)$ denotes a smoothly varying unit normal vector field along
$\gamma$, then the \emph{normal repulsion} $\rho$ of $\gamma$ over
the time interval $[t_{0},t]$ is given by \cite{var_theory} 
\begin{equation}
\rho(r,n)=\frac{1}{\sqrt{\left<n,C^{-1}(r)n\right>}},\label{eq:rhodef}
\end{equation}
measuring at time $t$ the normal component of the linearly advected
normal vector $\nabla F_{t_{0}}^{t}(r)n$ (see Fig. \ref{fig:rho}).
If $\rho>1$ pointwise along $\gamma$, then the the evolving material
line $F_{t_{0}}^{t}(\gamma)$ is repelling. Similarly, if $\rho<1$
holds pointwise along $\gamma$, then the the evolving material line
$F_{t_{0}}^{t}(\gamma)$ is attracting. 

Hyperbolic Lagrangian coherent structures (LCSs) are pointwise most
repelling or most attracting material lines with respect to small
perturbations to their tangent spaces \cite{var_theory,farazmand11,stretchlines}. Repelling and attracting
LCSs, respectively, are obtained as special trajectories of the differential
equations 
\begin{equation}
\dot{r}=\xi_{1}(r),\qquad\dot{r}=\xi_{2}(r),\label{eq:strain_stretch}
\end{equation}
that stay bounded away from points where $\xi_{i}$ cease to be well-defined.
These degenerate points $x_{0}$ are singularities of the Cauchy--Green
tensor field, satisfying $C(x_{0})=I$. The trajectories of the differential
equations in \eqref{eq:strain_stretch} are called strainlines and
stretchlines, respectively \cite{computeVariLCS,stretchlines}.
Computing the definition of $\rho$ in \eqref{eq:rhodef}, we obtain
that strainlines repel at a local rate of $\rho(r,n)=\sqrt{\lambda_{2}(r)}$,
and stretchlines attract at a rate of $\rho(r,n)=\sqrt{\lambda_{1}(r)}$.
Following the terminology used in the scientific visualization community
\cite{Ten_lines_2D,Tricoche}, we will refer to strainlines and stretchlines
collectively as \emph{tensorlines.}

A pointwise measure of how close a material curve is to being neutrally
stable is the \emph{neutrality} $\mathcal{N}(r,n)$, defined as 
\begin{equation}
\mathcal{N}(r,n)=\left(\rho(r,n)-1\right)^{2}.\label{eq:neutrality}
\end{equation}
Given the explicit normals known for tensorlines, their neutrality
can be computed as a sole function of the location $r$, and can be
written as 
\[
\mathcal{N}_{\xi_{1}}(r)=\left(\sqrt{\lambda_{2}(r)}-1\right)^{2},\qquad\mathcal{N}_{\xi_{2}}(r)=\left(\sqrt{\lambda_{1}(r)}-1\right)^{2},
\]
respectively, for strainlines and stretchlines.

In this paper, we will be seeking generalized non-twist curves (or
jet-cores) that are as close to neutral ($\mathcal{N}\equiv0$) as
possible. Requiring strictly zero neutrality along a material curve
$\gamma$ would, however, lead to an overdetermined problem. Indeed,
a material line with neutral stability at all its points would be
non-generic in an unsteady flow. Instead, we will be interested in
material lines that are close to minimizing the neutrality, while
also satisfying a minimal-shearing principle to be discussed later. 

Here we only work out a close-to-neutral condition for tensorlines,
as they will turn out to have special significance in our search for
shearless barriers. First, we define the convexity sets $\mathcal{C}_{\xi_{i}}$
of strainlines and stretchlines, respectively, as
\begin{equation*}
\mathcal{C}_{\xi_{i}}=\left\{ x_{0}\in U: \,\left<\partial_{r}^{2}\mathcal{N}_{\xi_{i}}(x_{0})\xi_{j}(x_{0}),\xi_{j}(x_{0})\right>>0,\qquad i\neq j\right\} ,\qquad i=1,2.
\end{equation*}

These sets are simply composed of points where the corresponding
neutrality is a convex function. We say that a compact tensorline segment
$\gamma$ \emph{is a weak minimizer} of its corresponding neutrality
$\mathcal{N}_{\xi_{i}}(r)$ if both $\gamma$ and the nearest trench
of $\mathcal{N}_{\xi_{i}}(r)$ lie in the same connected component
of $\mathcal{C}_{\xi_{i}}$. More specifically, a weak minimizer $\gamma$
of $\mathcal{N}_{\xi_{i}},$ with parametrization $r_{0}(s)$ and
smooth unit normal vector field $n_{0}(s)$, satisfies the condition
\begin{equation}
r_{0}(s)+\epsilon n_{0}(s)\in\mathcal{C}_{\xi_{i}},\qquad s\in[0,\sigma],\quad\epsilon\in[0,\epsilon_{0}(s)],\label{eq:convexity}
\end{equation}
where
\[
\epsilon_{0}(s)=arg\min\left\{ \left|\epsilon\right|\in(0,\infty)\,:\,\partial_{\epsilon}\mathcal{N}_{\xi_{i}}\left(r_{0}(s)+\epsilon n_{0}(s)\right)=0,\quad\partial_{\epsilon}^{2}\mathcal{N}_{\xi_{i}}\left(r_{0}(s)+\epsilon n_{0}(s)\right)>0\right\} .
\]

\section{Eulerian and Lagrangian shear}

For the steady two-dimensional steady flow shown in Fig. \ref{fig:canShearJet_tracers},
the classic Eulerian shear in the $x$ direction is defined as the
derivative of the horizontal velocity field in the vertical direction,
i.e., 
\begin{equation}
\frac{\partial u}{\partial y}=-2y,\label{Eul_shear}
\end{equation}
which vanishes on the line $y_{0}=0$. This line plays the role of
a jet core with a distinguished impact on tracer blobs in comparison
to other horizontal streamlines (see Figure \ref{fig:canShearJet_tracers}). 

The Eulerian shear, as the normal derivative of a velocity component
of interest, can certainly be computed for unsteady flows as well,
and is indeed broadly used in fluid mechanics \cite{batchelor}. However, instantaneously
shearless curves no longer act as invariant manifolds in the flow,
and thus will generally not create the characteristic tracer patterns
seen in Fig. \ref{fig:canShearJet_tracers}. As a result, the mathematical
description and systematic extraction of jet-core type material barriers
in unsteady flows has been an open problem, despite their ubiquitous
presence in plasma and geophysics.

To set the stage for a general description of jet-core-type structures,
we first need a Lagrangian definition of shear that captures the type
of material evolution seen in Fig. \ref{fig:canShearJet_tracers}
even in an unsteady flow. For an arbitrary material curve $\gamma$,
we select a parametrization $r(s)$ with $s\in[0,\sigma]$ for $\gamma$
at time $t_{0}$, and with the tangent vectors denoted as $r^{\prime}(s)$. 

We denote by $p(s)$ the pointwise tangential shear experienced over
the time interval $[t_{0},t]$ along the trajectory starting at time
$t_{0}$ from the point $r(s)$. Following \cite{geo_theory}, we
define this tangential shear in the Lagrangian frame as the $\gamma$-tangential
projection that a unit vector $n(s)=\left[r^{\prime}(s)\right]^{\perp}$
initially normal to $\gamma$ at $r(s)$ develops by time $t$, as
it is advected forward by the linearized flow $\nabla F_{t_{0}}^{t}(r(s))$
(see Fig. \ref{fig:rho}). Specifically, the \emph{Lagrangian shear} $p(s)$
is given by 
\begin{align}
p(s) & =\left\langle \frac{\nabla F_{t_{0}}^{t}(r(s))r^{\prime}(s)}{\left|\nabla F_{t_{0}}^{t}(r(s))r^{\prime}(s)\right|},\nabla F_{t_{0}}^{t}(r(s))\frac{\left[r^{\prime}(s)\right]^{\perp}}{\left|\left[r^{\prime}(s)\right]^{\perp}\right|}\right\rangle \nonumber \\
 & =\frac{\left\langle r^{\prime}(s),D(r(s))r^{\prime}(s)\right\rangle }{\sqrt{\left\langle r^{\prime}(s),C(r(s))r^{\prime}(s)\right\rangle \left\langle r^{\prime}(s),r^{\prime}(s)\right\rangle }},\label{eq:pdef}
\end{align}
where $\langle\cdot,\cdot\rangle$ denotes the Euclidean inner product,
and the tensor field $D$ is defined as 
\begin{equation}
D(x_{0})=\frac{1}{2}[C(x_{0})\Omega-\Omega C(x_{0})],\qquad\Omega=\left(\begin{array}{cc}
0 & -1\\
1 & 0
\end{array}\right).\label{eq:Ddef}
\end{equation}

\section{Variational principle for shearless transport barriers}

We seek generalized shearless curves as centerpieces of regions with
no observable variability in the averaged material shear. Assume that
$\epsilon>0$ is a minimal threshold above which we can physically
observe differences in material shear over the time interval $[t_{0},t].$
By smooth dependence on initial fluid positions, we will typically
observe an $\mathcal{O}(\epsilon)$ variability in shear within an
$\mathcal{O}(\epsilon)$-thick strip around a randomly chosen material
curve $\gamma$.  Our interest, however, is exceptional $\gamma$
curves around which $\mathcal{O}(\epsilon)$-thick coherent strips
show no observable variability in their average shearing. 

Based on the definition \eqref{eq:pdef}, the averaged Lagrangian
shear experienced along $\gamma$ over the time interval $[t_{0},t]$
can be written as
\begin{equation}
P(\gamma)=\frac{1}{\sigma}\int_{0}^{\sigma}p(s)\,\mathrm{d}s.
\label{eq:avLagShear}
\end{equation}

As we argued above, if an observable non-shearing material strip exists
around $\gamma$, then on $\epsilon$-close material curves we must
have $P(\gamma+\epsilon h)=P(\gamma)+\mathcal{O}(\epsilon^{2}),$
where $\epsilon h(s)$ denotes a small perturbation to $r(s)$. This
is only possible if the first variation of $P$ vanishes on $\gamma$:
\begin{equation}
\delta P(\gamma)=0.\label{eq:zerovari}
\end{equation}

This condition leads to the following weak form of the Euler-Lagrange
equation:
\begin{equation}
\delta P(\gamma)=\left[\left\langle \partial_{r^{\prime}}p,h\right\rangle \right]_{0}^{\sigma}+\int_{0}^{\sigma}\left[\partial_{r}p-\frac{d}{ds}\partial_{r^{\prime}}p\right]h(s)\, ds=0.\label{eq:vari}
\end{equation}

\section{Boundary conditions\label{sec:Boundary-conditions}}

We are interested in two types of boundary conditions for the
variational problem \eqref{eq:vari}:

\subsection{Variable endpoint boundary conditions}

Variable endpoint boundary conditions mean that $\gamma$
is a stationary curve with respect to all admissible perturbations,
i.e., it is the most observable type of centerpiece for shearless
coherent strip. As we show in Appendix \ref{bdrycond}, the only possible locations for variable
endpoint boundary conditions are those satisfying
\begin{equation}
C\left(r(0)\right)=C\left(r(\sigma)\right)=I.\label{eq:variBC}
\end{equation}

For completeness, we also consider another variable boundary condition in Appendix \ref{bdrycond} which results in non-zero Lagrangian shear \eqref{eq:pdef} and hence are not discussed here. 

\subsection{Fixed endpoint boundary conditions}

Fixed endpoint boundary conditions mean that $\gamma$ is
a stationary curve with respect to all perturbations that leave its
endpoints fixed. In this case, we have 
\begin{equation}
h(0)=h(\sigma)=0.\label{eq:homBC}
\end{equation}
 These boundary conditions do not place restrictions on the admissible
endpoints of $\gamma$. At the same time, a stationary curve under
these boundary conditions is generally expected to be less robust
or prevalent as a transport barrier than its variable-endpoint counterparts,
given that it only prevails as a stationary curve under a smaller
class of perturbations.

\section{Equivalent geodesic formulation: hyperbolic and parabolic barriers }

Under the above two boundary conditions, we obtain from \eqref{eq:vari}
the classic strong form of the Euler--Lagrange equations: 
\begin{equation}
\partial_{r}p-\frac{d}{ds}\partial_{r^{\prime}}p=0,\label{eq:EL}
\end{equation}
a complicated second-order differential equation for $r(s)$ . 

As we show in Appendix \ref{app:B}, however, any $\gamma$ satisfying
\eqref{eq:EL} also satisfies 
\begin{equation}
\delta\mathcal{P}_{\mu}(\gamma)=0,\qquad\qquad\mathcal{P}_{\mu}(\gamma)=\int_{\gamma}H_{\mu}(r(s),r^{\prime}(s))\,\mathrm{d}s,\qquad H_{\mu}(r(s),r^{\prime}(s))\equiv0,\label{eq:qdef-1}
\end{equation}
and hence represents a zero-energy stationary curve for the shear-energy-type
functional 
\begin{align}
H_{\mu}(r,r^{\prime}) & =\langle r^{\prime},D(r)r^{\prime}\rangle-\mu\sqrt{\left\langle r^{\prime},C(r)r^{\prime}\right\rangle \left\langle r^{\prime},r^{\prime}\right\rangle }\label{eq:Elambda}\\
\end{align}
for some choice of the parameter $\mu$. 

Of special interest to us is the case of pointwise shearless curves,
which we call \emph{perfect shearless barriers}. Such barriers should
prevail as influential transport barriers at arbitrary small scales.
Using the definition of the Lagrangian shear in \eqref{eq:pdef},
we conclude that curves with pointwise zero shear within the $H_{\mu}(r(s),r^{\prime}(s))\equiv 0$
energy surface all correspond to the parameter value $\mu=0.$ 

For this value of $\mu$, zero-energy stationary curves of the functional
$\mathcal{P}_{0}(\gamma)$ are null-geodesics of the Lorentzian metric
\begin{equation}
g(u,v)=\langle u,D(x_{0})v\rangle,\label{eq:Dmetric}
\end{equation}
which has metric signature $(-,+)$ \cite{LorentzianGeom}. The metric $g$ vanishes on
its null-geodesics, and hence these null-geodesics satisfy the implicit
first-order differential equation 
\begin{equation}
\langle r^{\prime}(s),D(r(s))r^{\prime}(s)\rangle\equiv0.\label{eq:0geod}
\end{equation}
 A direct calculation shows that all solutions of\eqref{eq:0geod}
satisfy
\begin{equation}
r^{\prime}(s)\parallel\xi_{i}(r(s)),\qquad i=1\,\,\mathrm{or\,\,2,}\label{eq:shearless_vector}
\end{equation}
therefore we obtain the following result. 
\begin{thm}
Perfect shearless barriers are null-geodesics of the Lorentzian metric
$g$, which are in turn composed of tensorlines of the Cauchy--Green
strain tensor $C$.
\end{thm}

\subsection{Hyperbolic barriers}

The geodesic transport barrier theory developed in \cite{geo_theory}
proposed that hyperbolic LCS are individual strainlines and stretchlines
that are most closely shadowed by locally most compressing and stretching
geodesics, respectively, of the Cauchy--Green strain tensor $C$. 

By contrast, here we have obtained from our shearless variational
principle \eqref{eq:zerovari} that tensorlines of $C$ are null-geodesics
for the tensor $D$. Instead of comparing tensorlines to Cauchy--Green
geodesics, therefore, one may simply locate hyperbolic LCSs as null-geodesics
of $D$ that 
\begin{description}
\item [{H1}] stay bounded away from Cauchy--Green singularities (i.e.,
points where $C=I$), elliptic LCSs (see \cite{geo_theory}) and parabolic LCSs (see below).

\item [{H2}] admit an extremum for the averaged compression or stretching,
respectively, among all their neighbors. These averages can be computed
by averaging $\sqrt{\lambda_{1}(x_{0})}$ and $\sqrt{\lambda_{2}(x_{0})}$,
respectively, along strainlines and stretchlines. 
\end{description}
Condition (H1) is required to hold because material curves crossing
Cauchy--Green singularities points have zero tangential and normal
stretching rates at the singularities, and hence lose their strict
normal attraction or repulsion property. It implies that hyperbolic
LCSs must satisfy Dirichlet boundary conditions, and none of their
interior points can be Cauchy-Green singularities either. As a result,
individual hyperbolic LCS are expected to fall in the less robust
and prevalent class of shearless barriers, as discussed in Section
\ref{sec:Boundary-conditions}.

\begin{figure}[H]
\begin{center}
\includegraphics[width=0.6\textwidth]{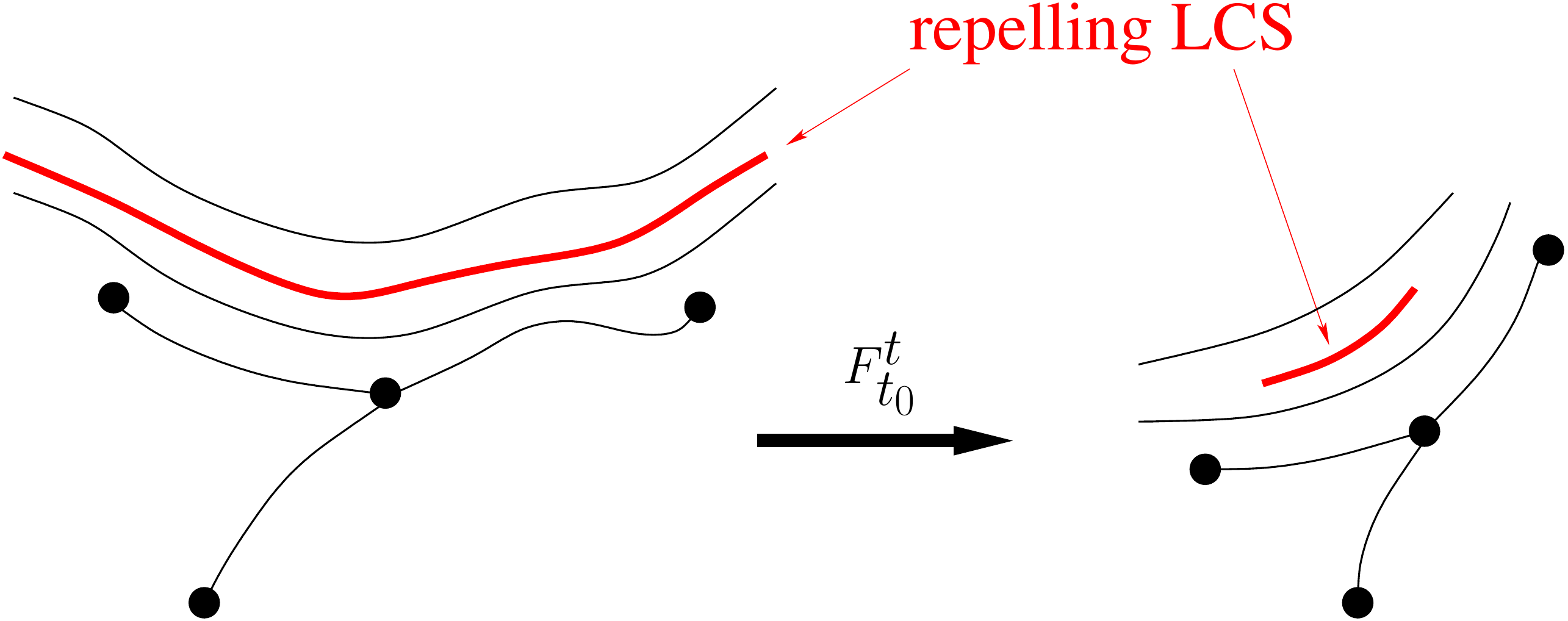}
\end{center}
\caption{Schematic representation of the properties of a repelling LCS (red) among nearby strainlines (black) and Cauchy--Green singularities (dots).  The repelling LCS stays away from singularities of Cauchy--Green singularities. While the length of any strainline shrinks as advected under the flow map $F_{t_0}^t$, the length of a repelling LCS shrinks more than any $C^1$-close strainline.}
\label{fig:hypLCS}
\end{figure}

Condition (H2) simply implements the definition of LCS as locally
most repelling or attracting material curves, reducing an originally
infinite-dimensional extremum problem to maximization within a one-dimensional
family of strainlines or stretchlines. We summarize the implications
of our shearless variational principle for hyperbolic LCS detection.
\begin{prop}
{[}Hyperbolic LCS as shearless barriers{]} Hyperbolic LCSs at time
$t_{0}$ are null-geodesics of the Lorentzian metric $g$ that are
bounded away from $C(x_{0})=I$ singularities of the Cauchy--Greens
strain tensor. In addition, repelling LCSs have an average stretching
smaller than that of any $C^{1}$ close null-geodesic of $g$ (see Fig. \ref{fig:hypLCS} for an illustration). 
Furthermore, attracting LCSs have an average stretching larger than that of any
$C^{1}$ close null-geodesic of $g$.
\end{prop}

\subsection{Parabolic barriers\label{sub:Parabolic-barriers}}

Our main focus is to find generalized jet cores in the Lagrangian
frame for unsteady flows of arbitrary time dependence. We shall refer
to such generalized jet cores here as \emph{parabolic transport barriers}.

The general solution \eqref{eq:shearless_vector} of our variational
principle certainly allows for further types of shearless barriers
beyond hyperbolic LCSs. These further barriers are also composed of
strainlines and stretchlines, but contain Cauchy--Green singularities
and hence fail to be hyperbolic material lines. As discussed in section
\eqref{sec:Boundary-conditions}, such non-hyperbolic barriers are
the most influential if they satisfy variable-endpoint boundary conditions
for our shearless variational principle, i.e., their endpoints are
Cauchy--Green singularities.

In addition, in order to provide a generalization of jet cores, we
are interested in non-hyperbolic shearless barriers that have no distinct
(repelling or attracting) stability type along their interior points.
To this end, we require parabolic barriers to be also weak minimizers
of their neutrality in the sense of Section \ref{sec:stability}. 

\begin{figure}[H]
\begin{center}
\includegraphics[width=0.3\textwidth]{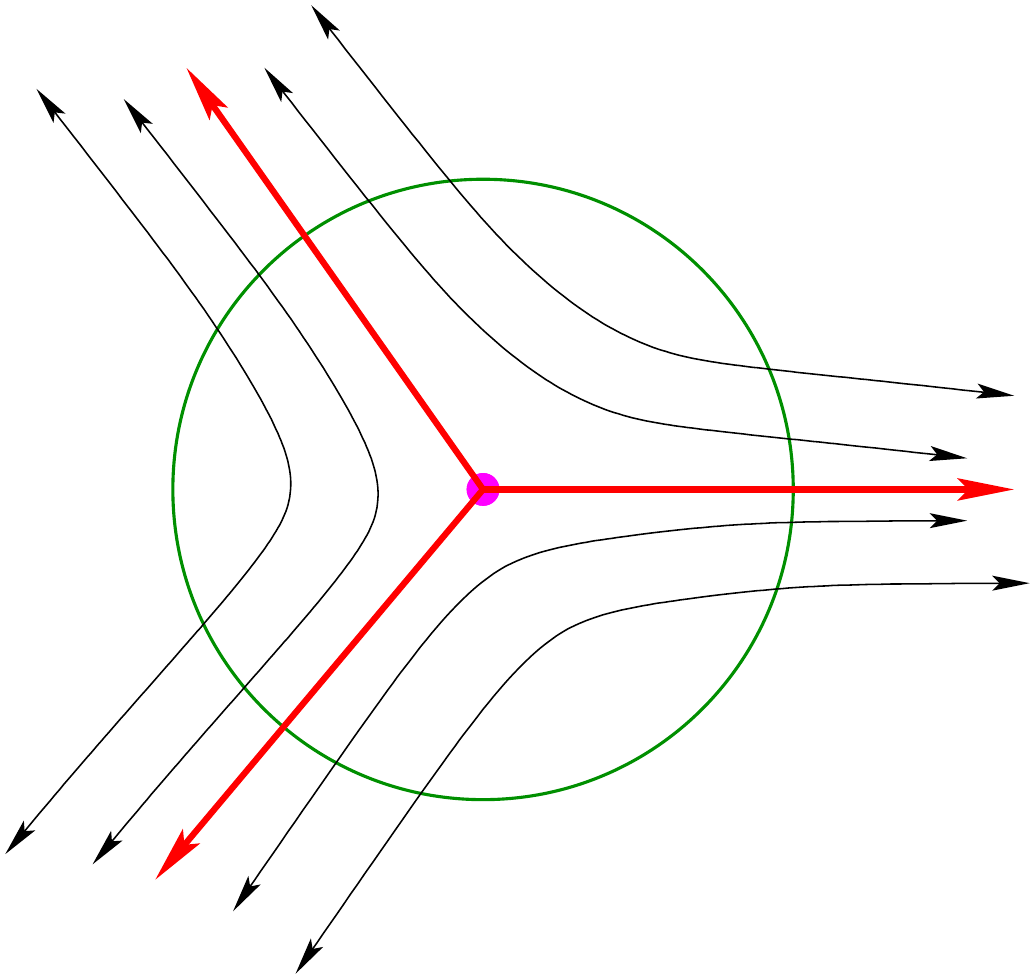}\hspace{.08\textwidth}
\includegraphics[width=0.25\textwidth]{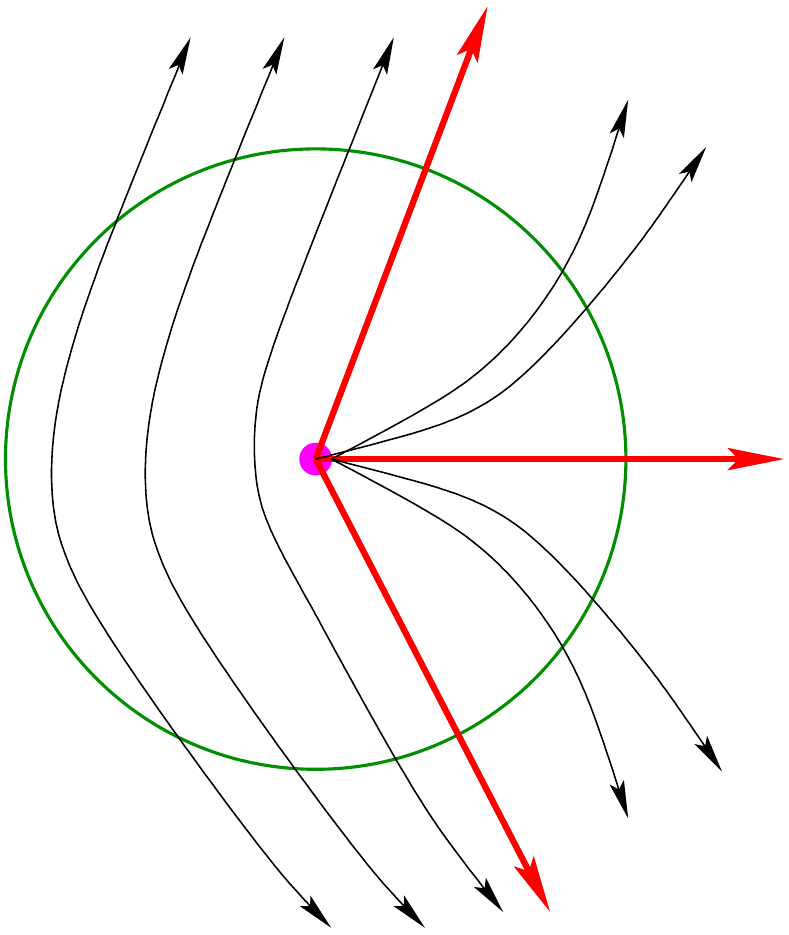}
\end{center}
\caption{Topology of tensorlines (black) around a trisector (left) and a wedge
(right) singularity (magenta). The tensorlines shown in red form the
separatrices.}
\label{fig:tri_wedge} 
\end{figure}

Finally, for reasons of physical relevance and observability, our
definition of a parabolic barrier will further restrict our consideration
to strainline--stretchline chains that are unique between the two
singularities they connect, and are also structurally stable with
respect to small perturbations. Based on our review of tensorline
singularities in Appendix \ref{app:C}, strainlines connecting singularities
are only structurally stable and unique if they connect a \emph{trisector}
singularity to a \emph{wedge} singularity (see Fig. \ref{fig:tri_wedge}).
An identical requirement holds for stretchlines.

\begin{figure}[H]
\begin{center}
\includegraphics[width=0.5\textwidth]{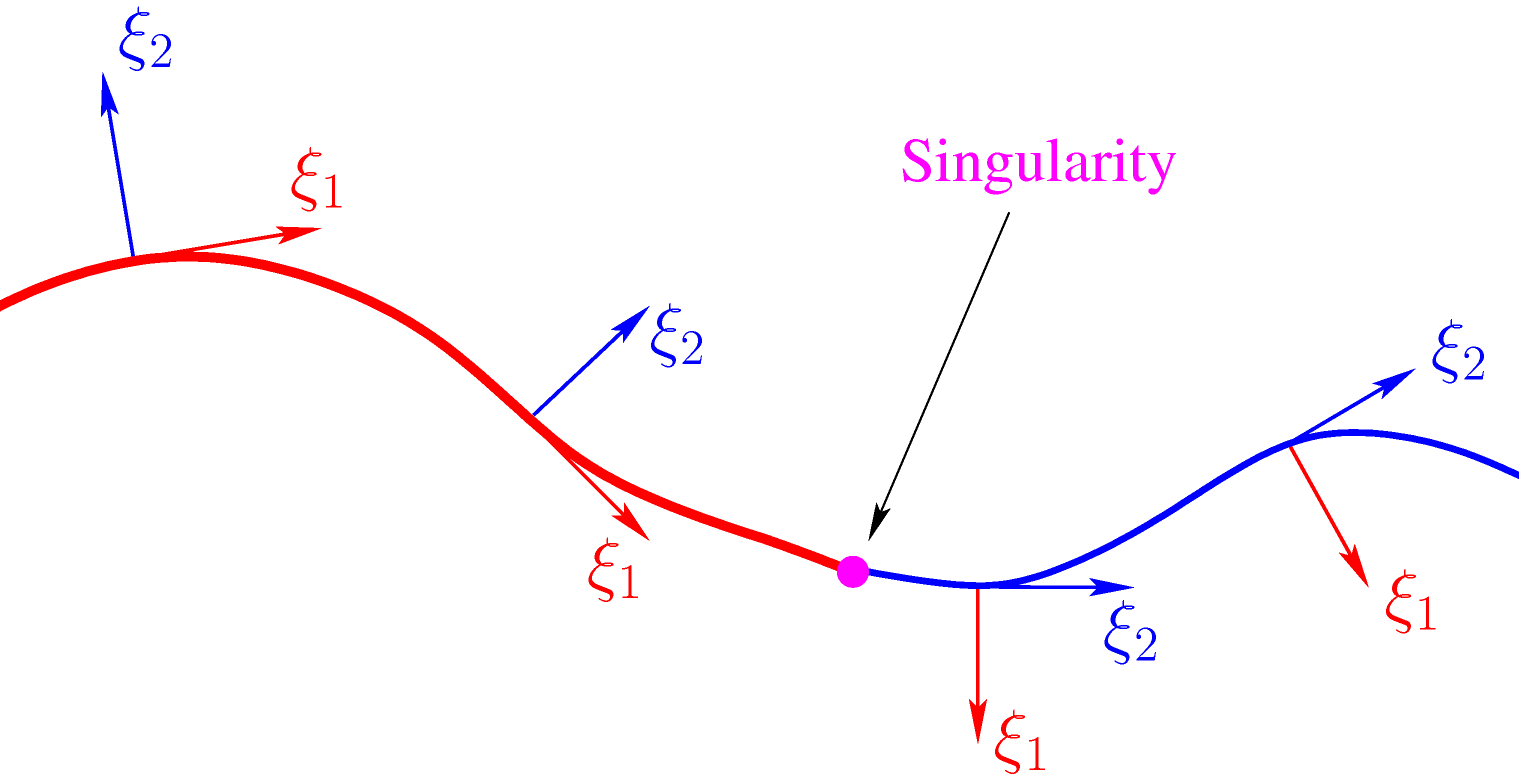}\\
\includegraphics[width=0.5\textwidth]{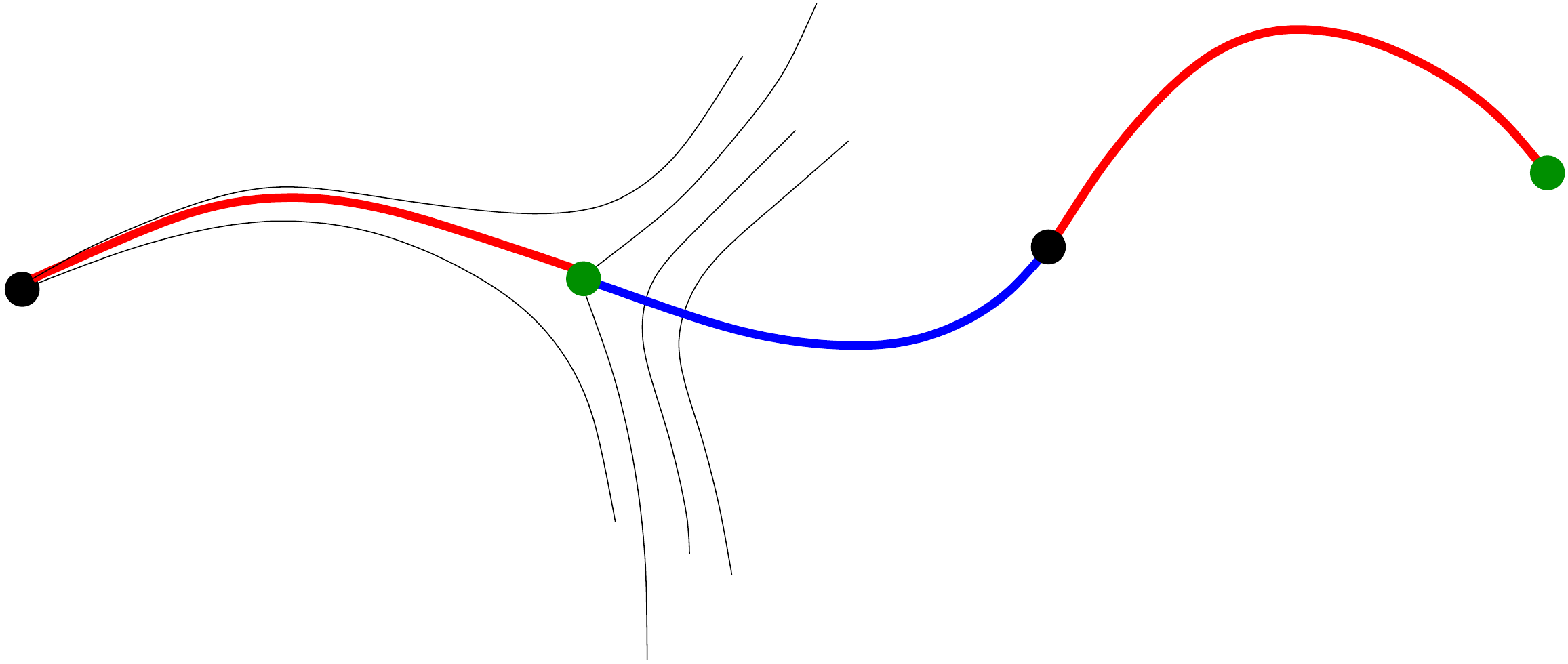}
\end{center}
\caption{Top: Smooth connection of strainlines (red curve) and stretchlines (blue curve) only occurs at Cauchy--Green singularities. 
Bottom: An alternating chain of strainlines (red) and stretchlines (blue) connecting trisectors (green) and wedges (black). An schematic phase portrait of strainlines (thin black lines) is shown around one of the trisector singularities. The strainline marked by red color is the unique connection between that trisector and the wedge on its left.}
\label{fig:tensorlines_chain}
\end{figure}

We then have the following definition.
\begin{defn}
{[}Parabolic barriers{]}\label{def:parabolic barriers} Let $\gamma$
denote the time $t_{0}$ position of a compact material line. Then
this material line is a \emph{parabolic transport barrier} over the
time interval $[t_{0},t]$ if the following two conditions are satisfied:\end{defn}
\begin{description}
\item [{P1}] $\gamma$ is an alternating chain of strainlines and stretchlines,
which is a unique connection between a wedge- and and a trisector-type
singularity of the tensor field $C(x_{0})$ (see Fig. \ref{fig:tensorlines_chain}).
\item [{P2}] Each strainline and stretchline segment in $\gamma$ is a
weak minimizer of its associated neutrality.
\end{description}

\begin{example}\label{ex:FTLE-Trench}
{[}\emph{An FTLE trench is not necessarily a parabolic barrier}{]}
Since our notion of a parabolic barrier requires a minimality condition
on $\lambda_{2}$, one may speculate whether a trench of the Finite-Time
Lyapunov Exponent (FTLE) field will always be a shearless barrier.
Such an approach of detecting jet cores by trenches of the combined
forward and backward FTLE field was considered in \cite{Javier_shear}.
While the trench of the FTLE field can indeed be an indicator of a
jet core, the following example of a steady two-dimensional incompressible
flow shows that this is not necessarily the case. Consider the incompressible
flow 
\begin{equation}
\begin{split} & \dot{x}=x\left(1+3y^{2}\right),\\
 & \dot{y}=-y-y^{3}.
\end{split}
\label{eq:ftle_counterEx}
\end{equation}

The line $y=0$ is an invariant, attracting set, yet numerical simulations
show that it is also a trench of the FTLE field, as seen in Figure
\ref{fig:ftle_counterEx}. The figure also shows
by tracer advection that this trench is a hyperbolic (attracting)
LCS, as opposed to a parabolic barrier acting as a jet core. 
\begin{figure}[H]
\begin{center}
\includegraphics[width=0.35\textwidth]{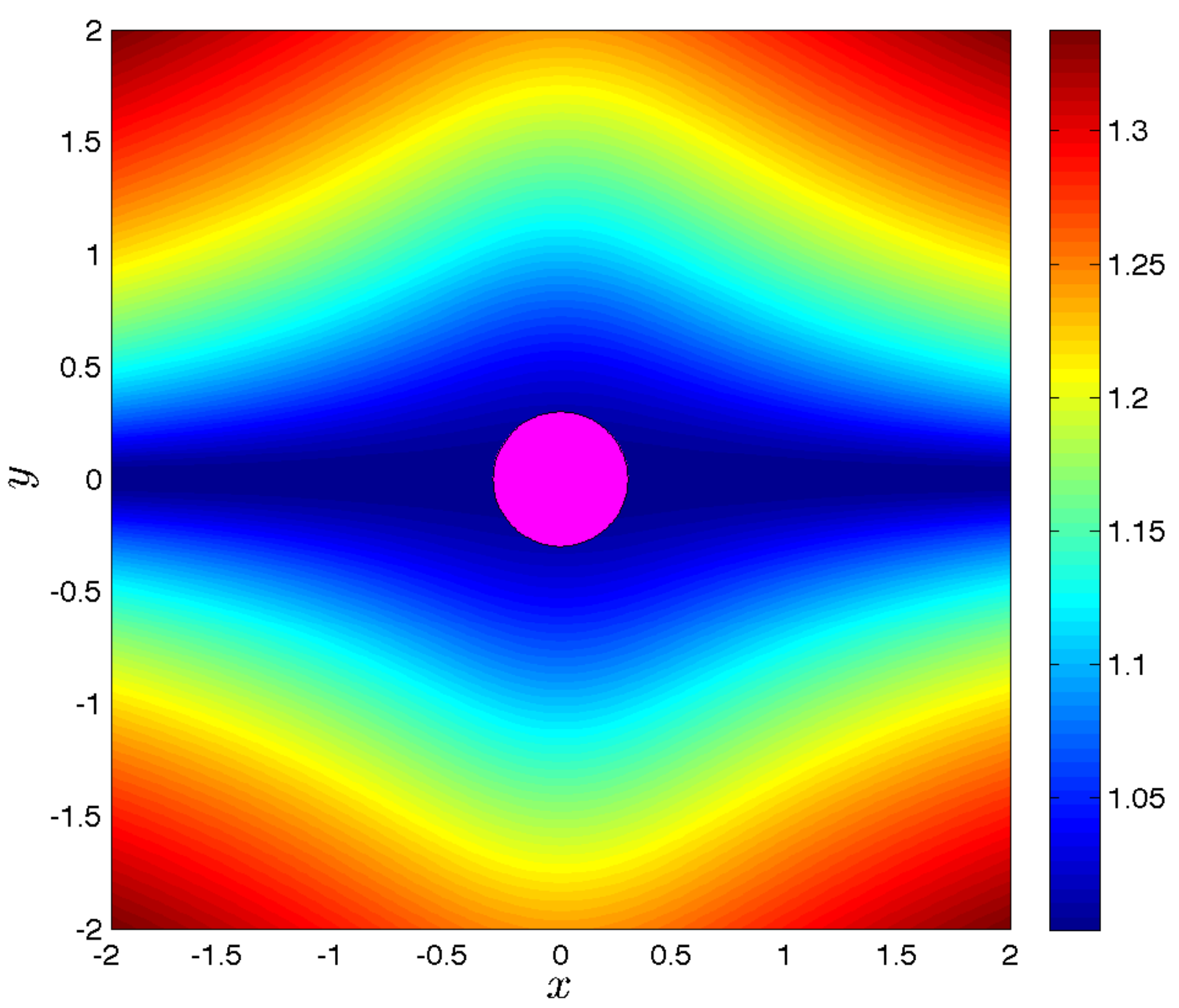}
\includegraphics[width=0.35\textwidth]{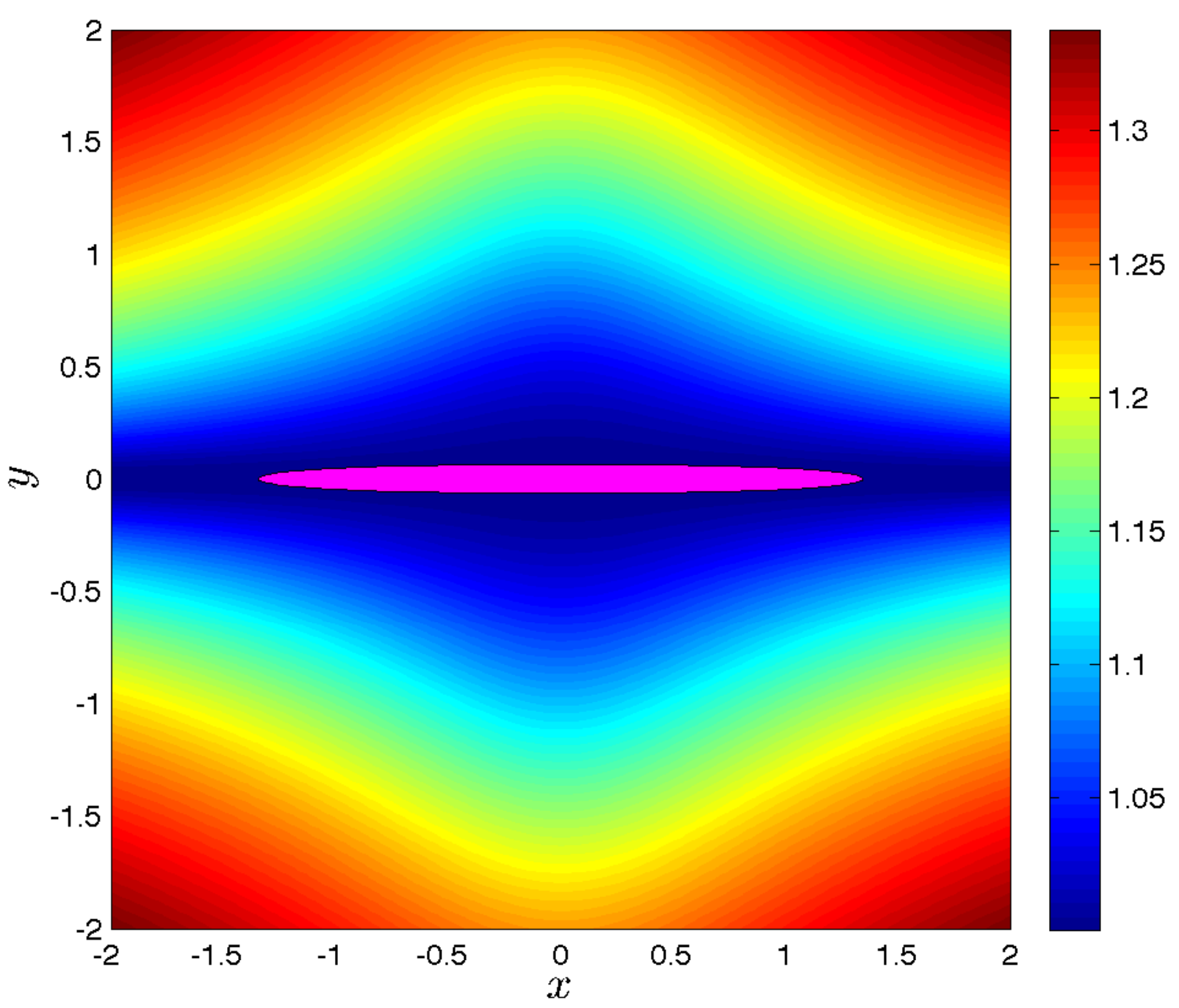}
\end{center}
\caption{The tracer evolution for system (\ref{eq:ftle_counterEx}). Left:
Initial circular blob of tracers centered at the origin at time $t=0$.
Right:  The advected tracer at time $t=1.5$. The forward-time FTLE
field with integration time $T=10$ is shown in the background. }
\label{fig:ftle_counterEx}
\end{figure}

\end{example}

\section{Automated numerical detection of parabolic barriers} \label{sec:methods} 
Definition \ref{def:parabolic barriers} provides the basis for the
identification of parabolic barriers in finite-time flow data. Using
the numerical details surveyed in Appendix \ref{app:C}, we implement conditions
P1 and P2 of Definition \ref{def:parabolic barriers} as follows:
\begin{enumerate}
\item Compute the Cauchy--Green strain tensor $C$ on a two-dimensional
grid in the $(x_{1},x_{2})$ variables.
\item Detect the singularities of $C$ by finding the common zeros of $f=C_{11}-C_{22}$
and $g=C_{12}$. 
\item For any trisector singularity of the $\xi_{1}$ vector field, follow
strainlines emanating from the singularity and identify among them
the separatrices connecting the trisectors to wedges. Repeat the same
procedure for the $\xi_{2}$ vector field to find trisector-wedge separatrices
among stretchlines.
\item Out of the computed separatrices, keep the strainline separatrices
satisfying $\left<\partial_{r}^{2}\mathcal{N}_{\xi_{1}}(x_{0})\xi_{2}(x_{0}),\xi_{2}(x_{0})\right>>0$,
and the stretchline separatrices satisfying $\left<\partial_{r}^{2}\mathcal{N}_{\xi_{2}}(x_{0})\xi_{1}(x_{0}),\xi_{1}(x_{0})\right>>0$. 
\item Build smoothly connecting, alternating stretchline-strainline heteroclinic
chains form the separatrices so obtained.
\item Finally, keep only the heteroclinic chains whose individual components
are weak minimizers of their neutralities.
\end{enumerate}

\section{Numerical examples}

\label{sec:results}

\subsection{Standard non-twist map}

\label{section:SNTM} We first consider the standard non-twist map
(SNTM) 
\begin{equation}
\begin{split} & x_{n+1}=x_{n}+a\left(1-y_{n+1}^{2}\right),\\
 & y_{n+1}=y_{n}-b\sin(2\pi x_{n}),
\end{split}
\label{eq:sntm}
\end{equation}
which was first studied in detail in \cite{Diego_NT}, and has since
become a generally helpful model in understanding shearless KAM curves
in two-dimensional steady or temporally periodic incompressible
flows. 

For $b=0$, the map \eqref{eq:sntm} is a discretized version of the
canonical parallel shear flow \eqref{fig:canShearJet_tracers} with
vanishing Eulerian shear along $y=0$. For steady perturbations of
\eqref{fig:canShearJet_tracers}, one still has a steady streamfunction
whose dynamics is integrable and the shearless barriers can be understood
as the lack of Hamiltonian twist. For $b\neq0$, the SNTM corresponds
to the evolution of a time-periodic perturbation of \eqref{fig:canShearJet_tracers}. 

\begin{figure}[H]
\begin{center}
\includegraphics[width=0.45\textwidth]{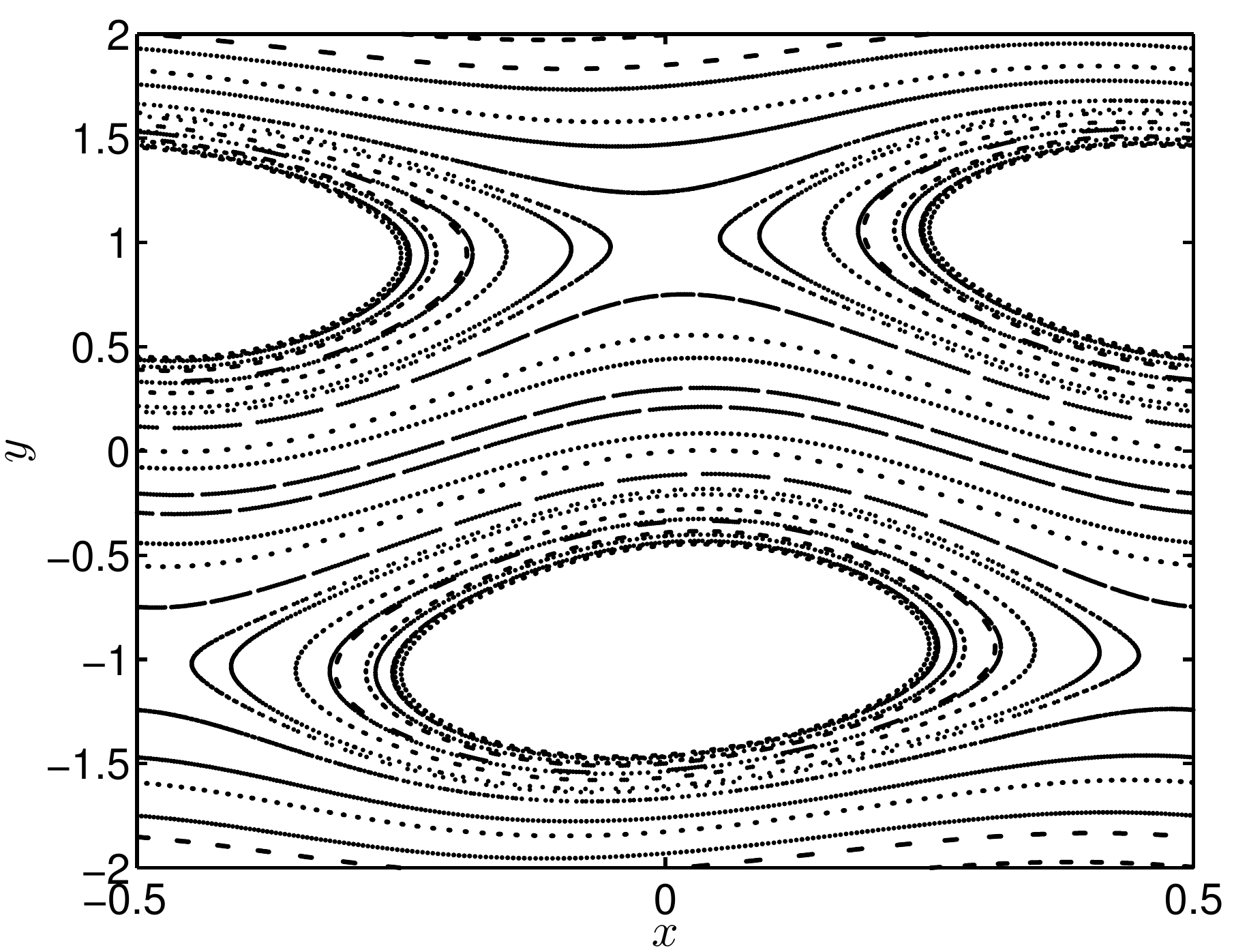} 
\includegraphics[width=0.45\textwidth]{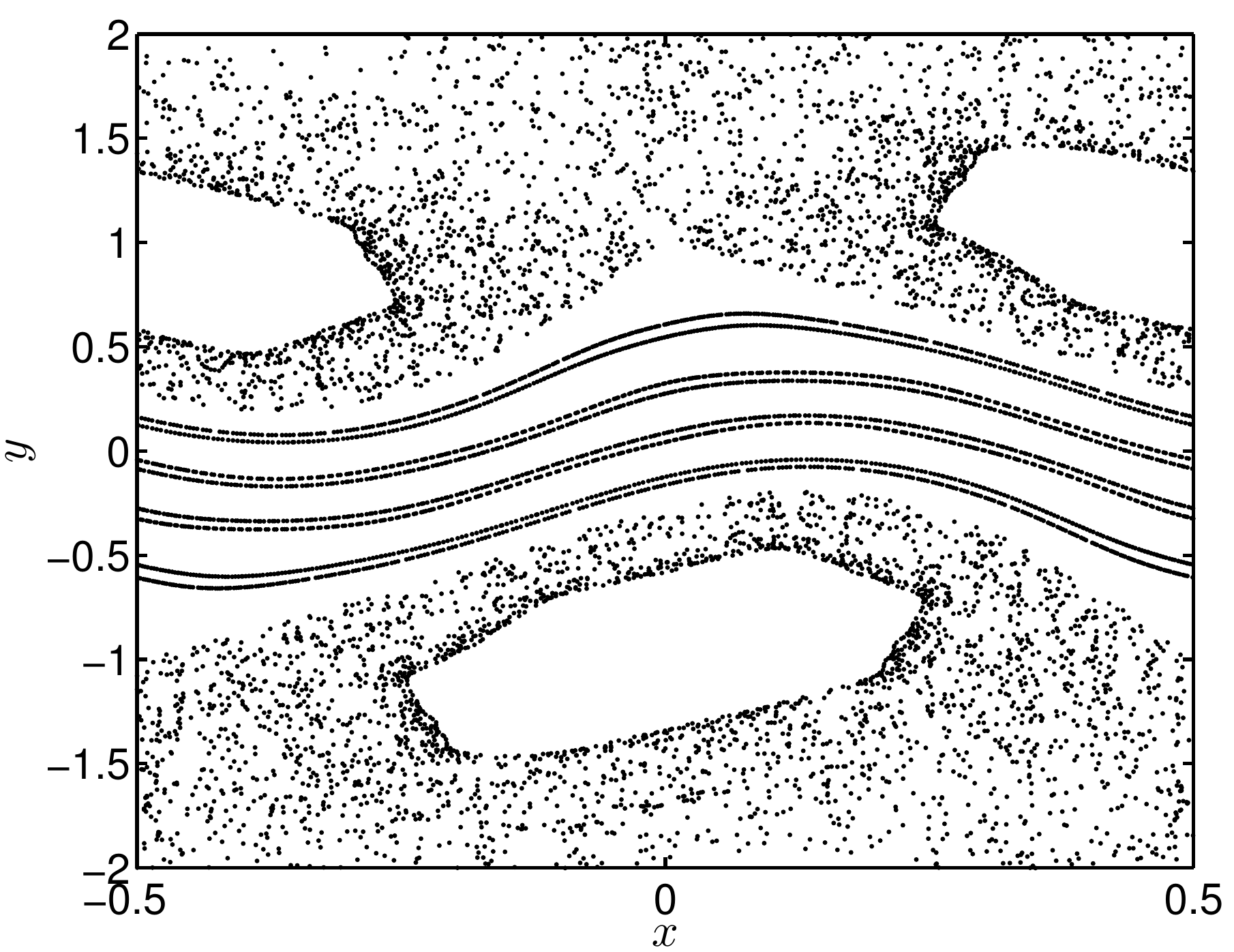}
\end{center}
\caption{The standard non-twist map. Left: Integrable parameters:
$a=0.08$, $b=0.125$ Right: Chaotic parameters: $a=0.27$, $b=0.38$.}
\label{fig:pmap_sntm} 
\end{figure}

For the parameter values $a=0.08$, $b=0.125$, the SNTM is integrable
and well-understood. We choose these parameters to illustrate the
performance of our theory and extraction methodology for parabolic
barriers. Figure \ref{fig:pmap_sntm}a shows the orbits of SNTM
for these integrable parameters.

In this integrable case, the location of shearless barriers is no
longer trivial, but can be found by the theory of indicator points \cite{Indicator_points}.
Specifically, initial conditions for the shearless barrier are given
by 
\begin{equation}
x=\pm\left(\frac{a}{2}+\frac{1}{4}\right)\ \ \mbox{and}\ \ y=0,\label{eq:ind_pnts}
\end{equation}
and the full barrier can be constructed by iterating these initial
condition under the map (\ref{eq:sntm}). Therefore, we can compare
the parabolic barrier computed from finitely many iterations using
the steps in Section \ref{sec:methods} with the exact asymptotic
shearless barrier of the map.

\begin{figure}[H]
\begin{center}
\includegraphics[width=0.35\textwidth]{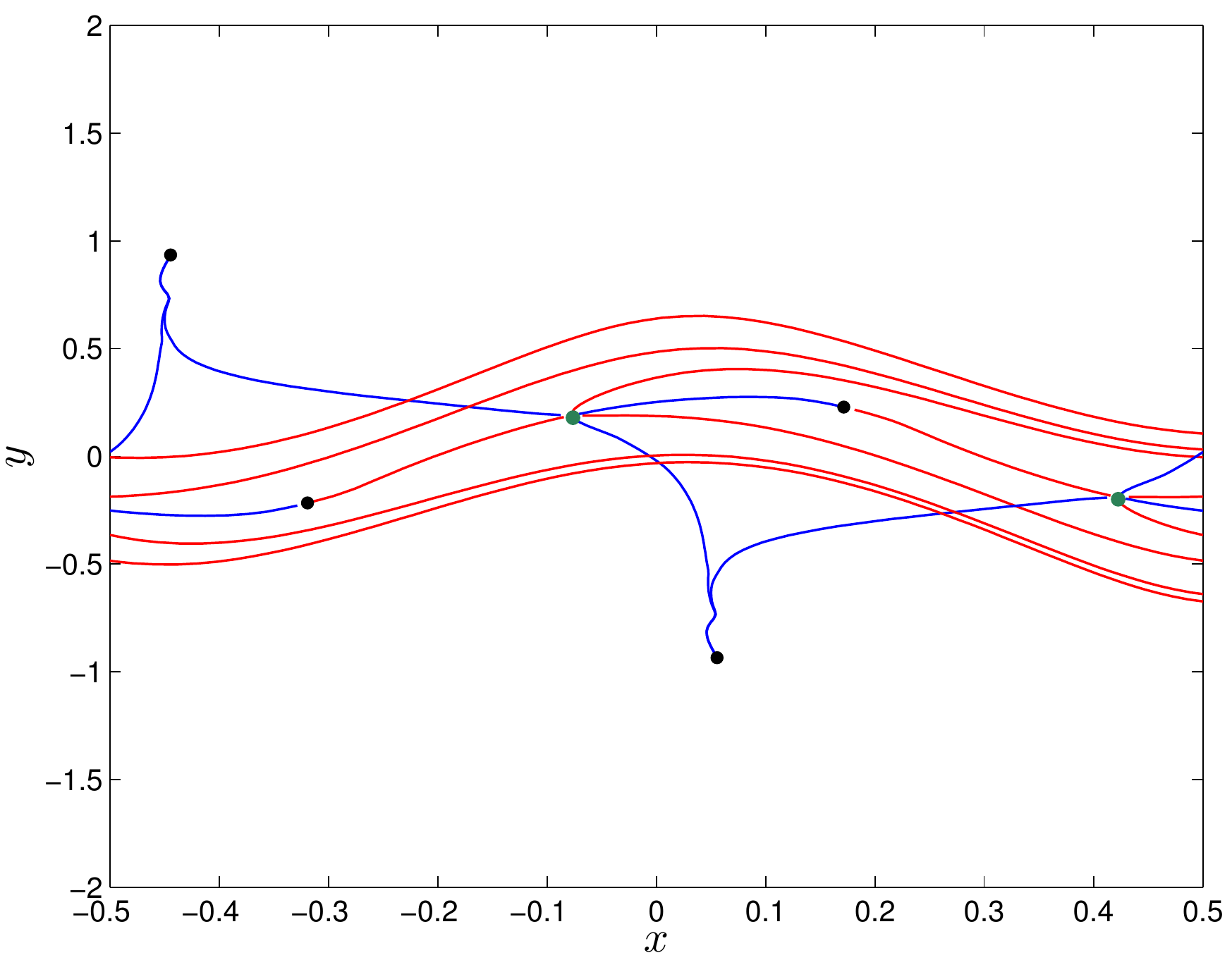}
\includegraphics[width=0.35\textwidth]{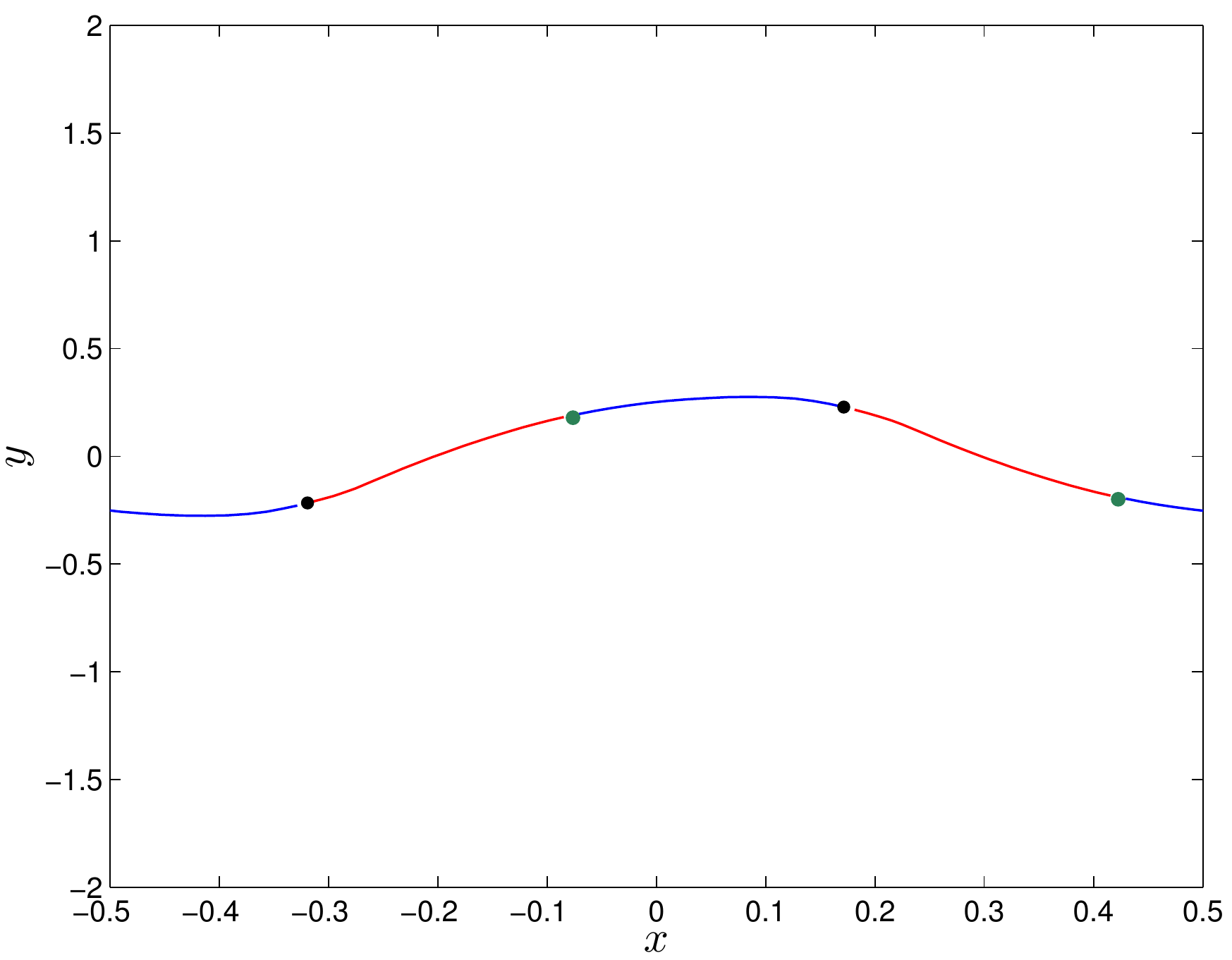}
\end{center}
\caption{Left: Heteroclinic tensorlines between trisector and wedge singularities
of the Cauchy--Green strain tensor in the integrable SNTM: strainlines
(red) and stretchlines (blue). The black and green dots mark the
wedge and trisector singularities, respectively. Right: The extracted
parabolic barrier consists of the single alternating sequence of tensorlines
that satisfy conditions P1-P2 of Definition \ref{def:parabolic barriers}.}
\label{fig:sntm_int_tensorlines} 
\end{figure}

Figure \ref{fig:sntm_int_tensorlines} shows all heteroclinic tensorlines
connecting trisectors to wedges (left panel). In the domain $[-0.5,0.5]\times[-2,2]$
and for $100$ iterations of the SNTM, we find $6$ singularities:
2 trisectors (green dots) and 4 wedges (black dots). Only 4 alternating
sequence of tensorlines satisfy conditions P1 and P2 of Definition
\ref{def:parabolic barriers}. Figure \ref{fig:sntm_int_tensorlines}
also shows the extracted parabolic barrier, i.e., a heteroclinic chain
formed by four tensorlines (note the periodicity in $x$). This parabolic
barrier represents the finite-time version of the exactly known asymptotic
shearless KAM curve.

\begin{figure}[H]
\begin{center}
\includegraphics[width=0.3\textwidth]{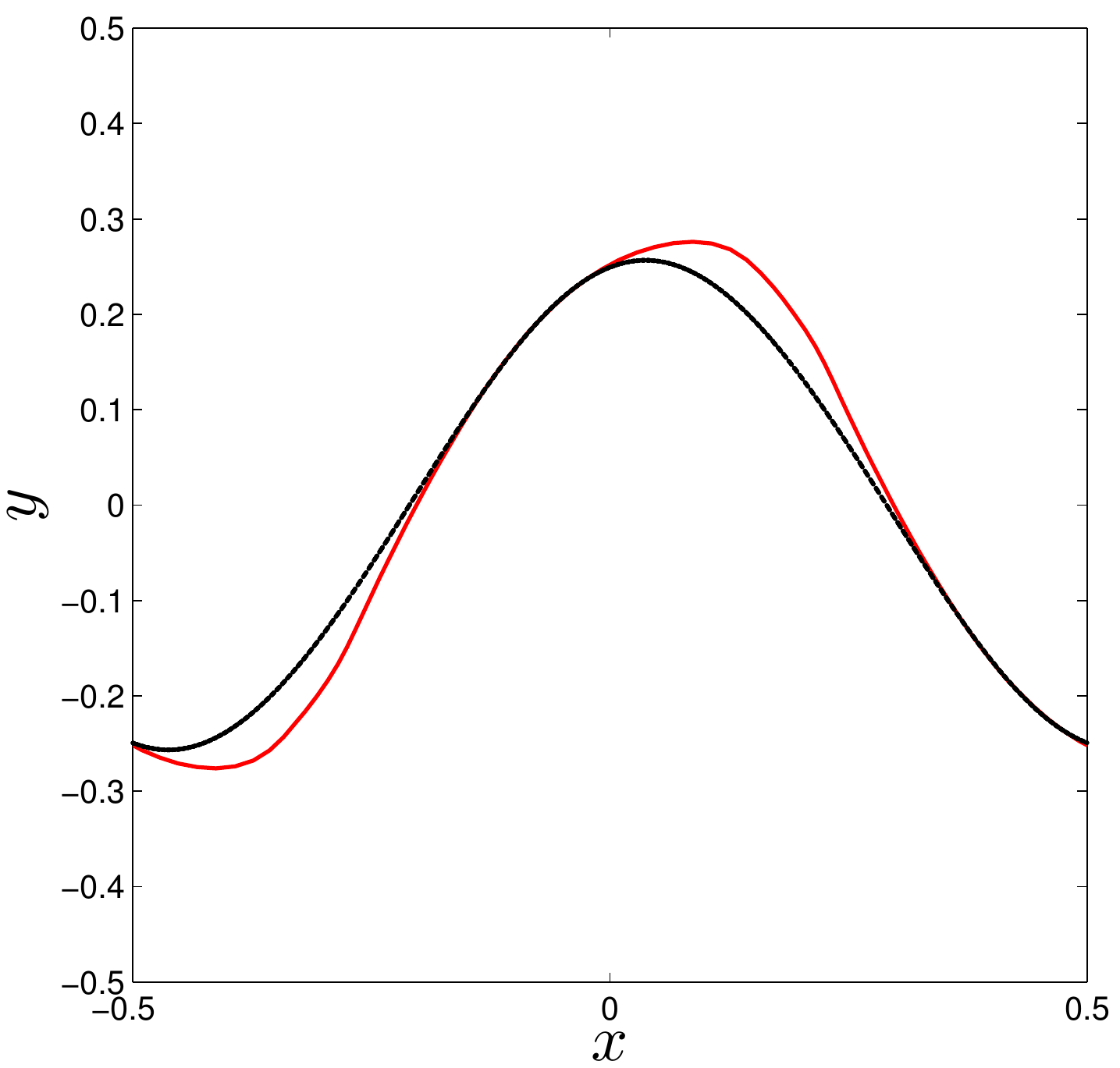}
\includegraphics[width=0.3\textwidth]{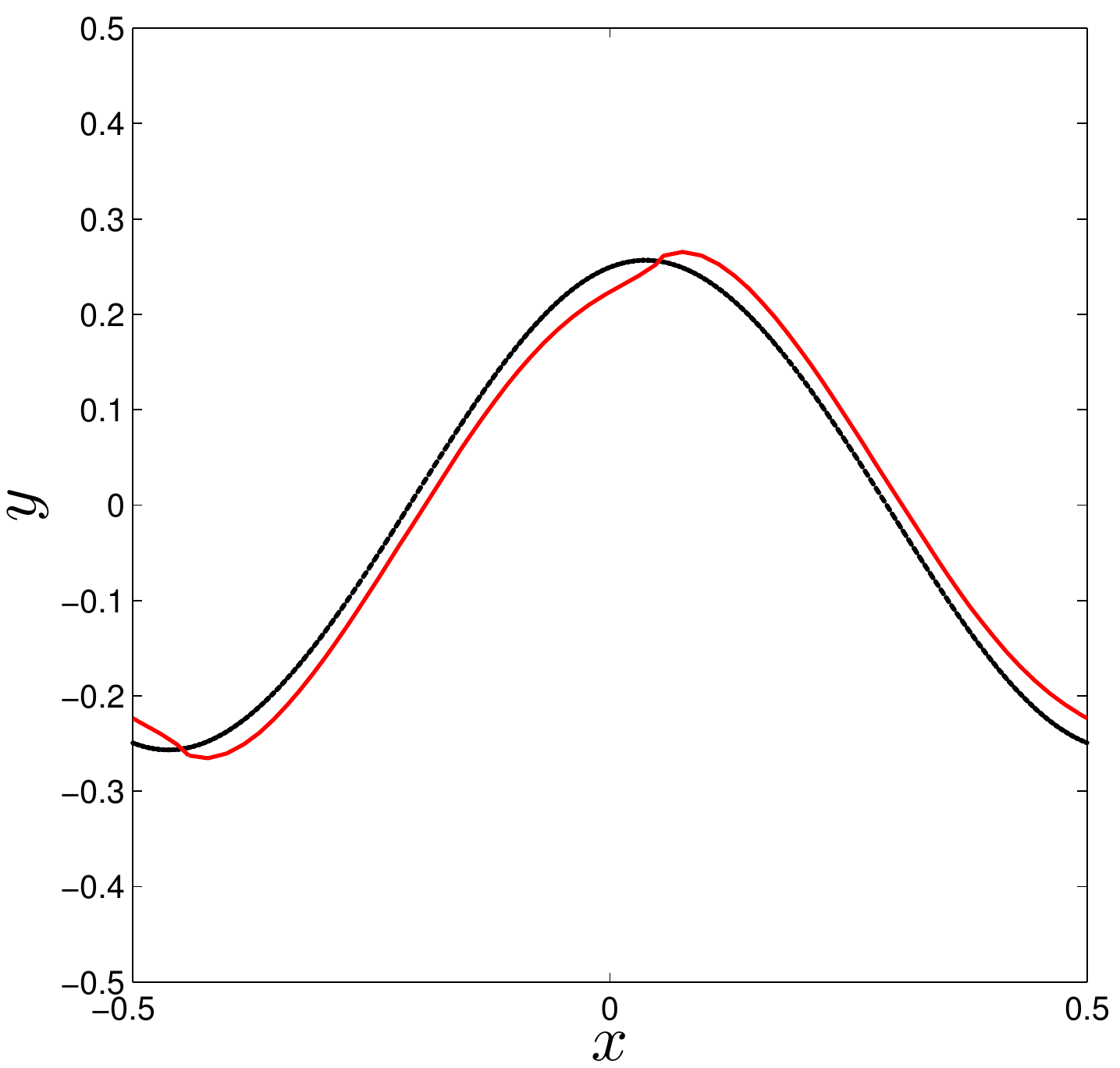}
\includegraphics[width=0.3\textwidth]{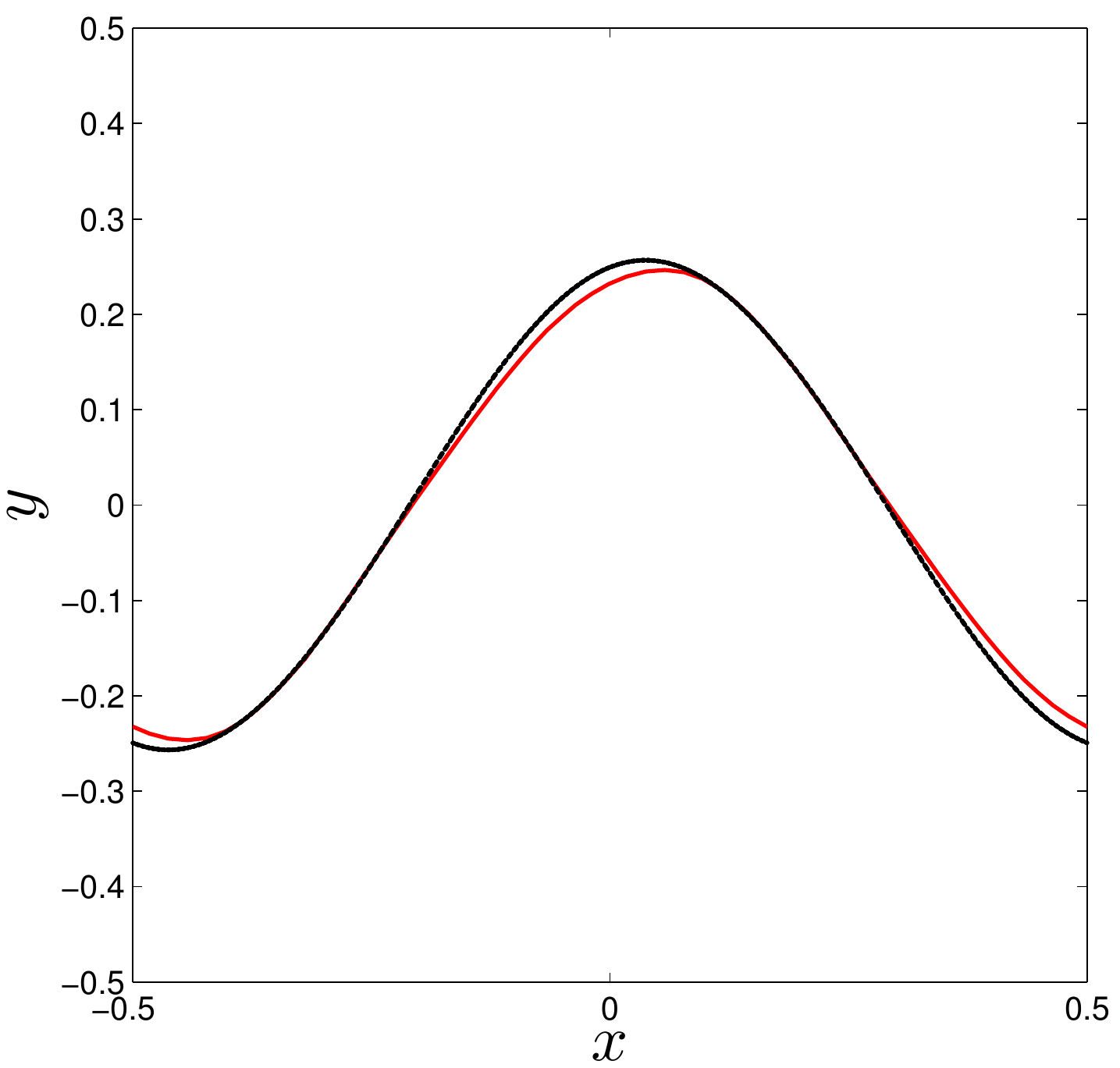}
\end{center}
\caption{The red curve shows the computed finite-time shearless barrier from
$100$ (left), $200$ (middle) and $300$ (right) iterations of the
integrable SNTM with parameters $a=0.08$ and $b=0.125$. The black
curve marks the exact location of the barrier.}
\label{fig:compare_sntm_int_aut} 
\end{figure}

One can also compute the parabolic barrier for higher iterations of
the SNTM map with the same procedure. As the number of iterations
increase, the computed parabolic barrier converges to the exact asymptotic
barrier. In Fig. \ref{fig:compare_sntm_int_aut},
we show this convergence up to $300$ iterations. For higher iterations,
the two barriers become practically indistinguishable. The exact barrier (black curve) 
in Fig. \ref{fig:compare_sntm_int_aut} is computed from $200$ iterations of the indicator points \eqref{eq:ind_pnts}.

\begin{figure}[H]
\begin{center}
\includegraphics[width=0.3\textwidth]{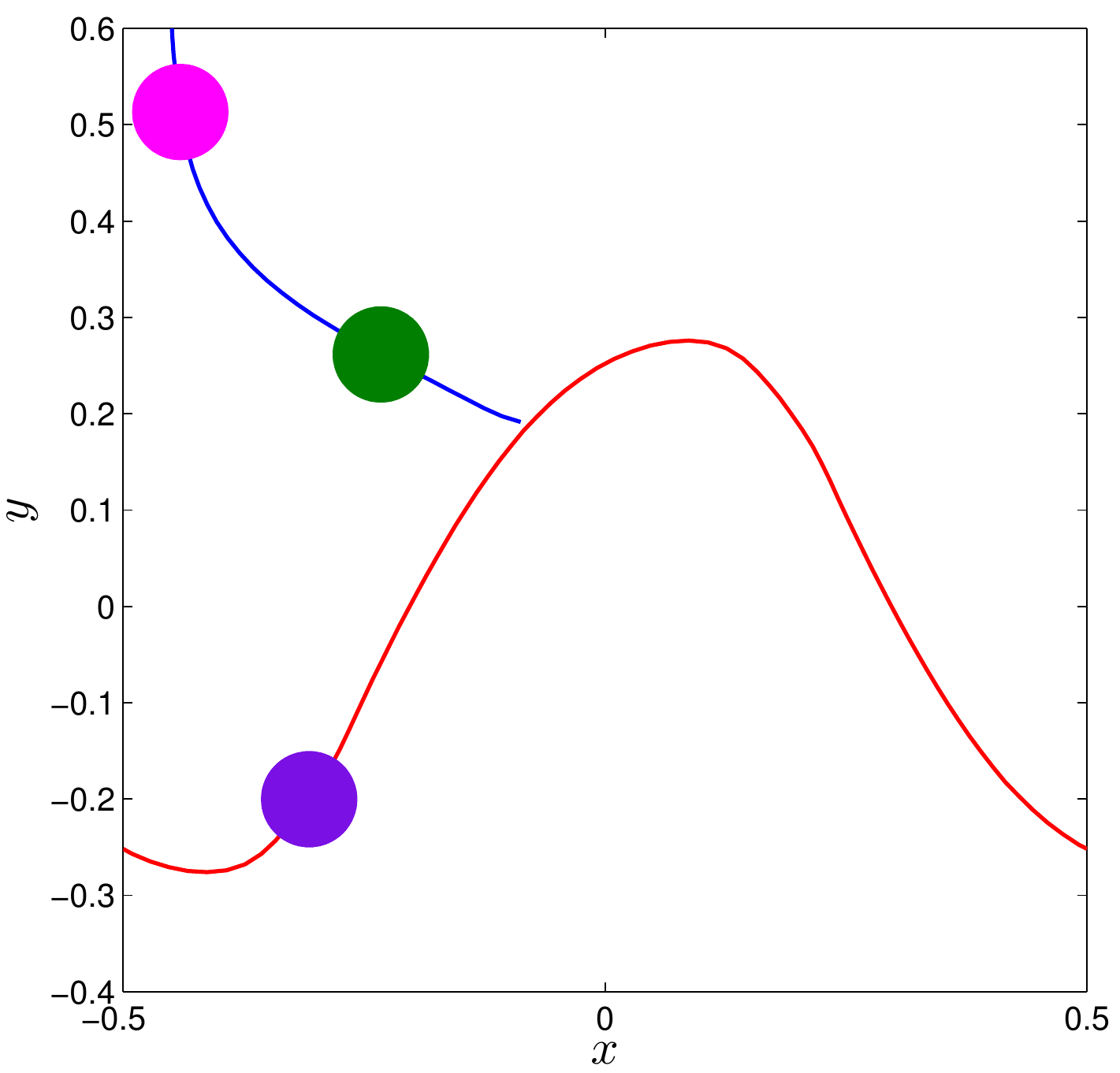}
\includegraphics[width=0.37\textwidth]{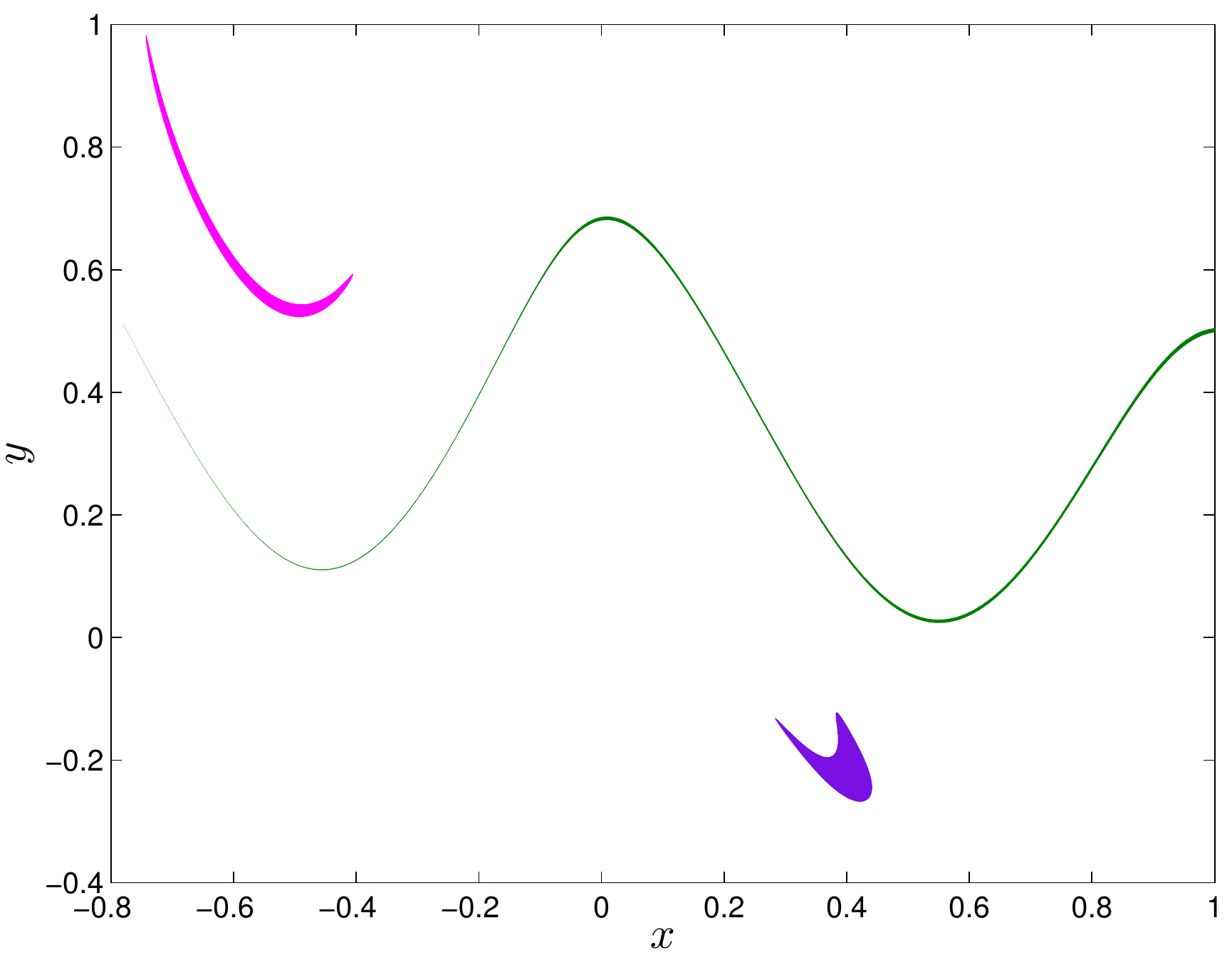}
\end{center}
\caption{Parabolic barrier and its impact on tracers in the integrable SNTM.}
\label{fig:tracers_sntm_int} 
\end{figure}

The evolution of circular tracers off and on the computed parabolic
barriers is shown in Fig. \ref{fig:tracers_sntm_int} The purple tracer
in the left plot of Fig. \ref{fig:tracers_sntm_int} is located on
the computed parabolic barrier (red). The magenta and green tracers
are centered on a tensorline (blue) that does not satisfy condition
P2 of Definition \ref{def:parabolic barriers}. The images of all
the tracers after $100$ iterations of SNTM are shown in the right panel
of Fig. \ref{fig:tracers_sntm_int}. While the purple
tracer undergoes a small boomerang-like deformation expected along
parabolic barriers (jet cores), the other two tracer blobs experience
substantial stretching. This illustrates that condition P2 is indeed
essential in identifying parabolic barriers.

\begin{figure}[H]
\begin{center}
\includegraphics[width=0.5\textwidth]{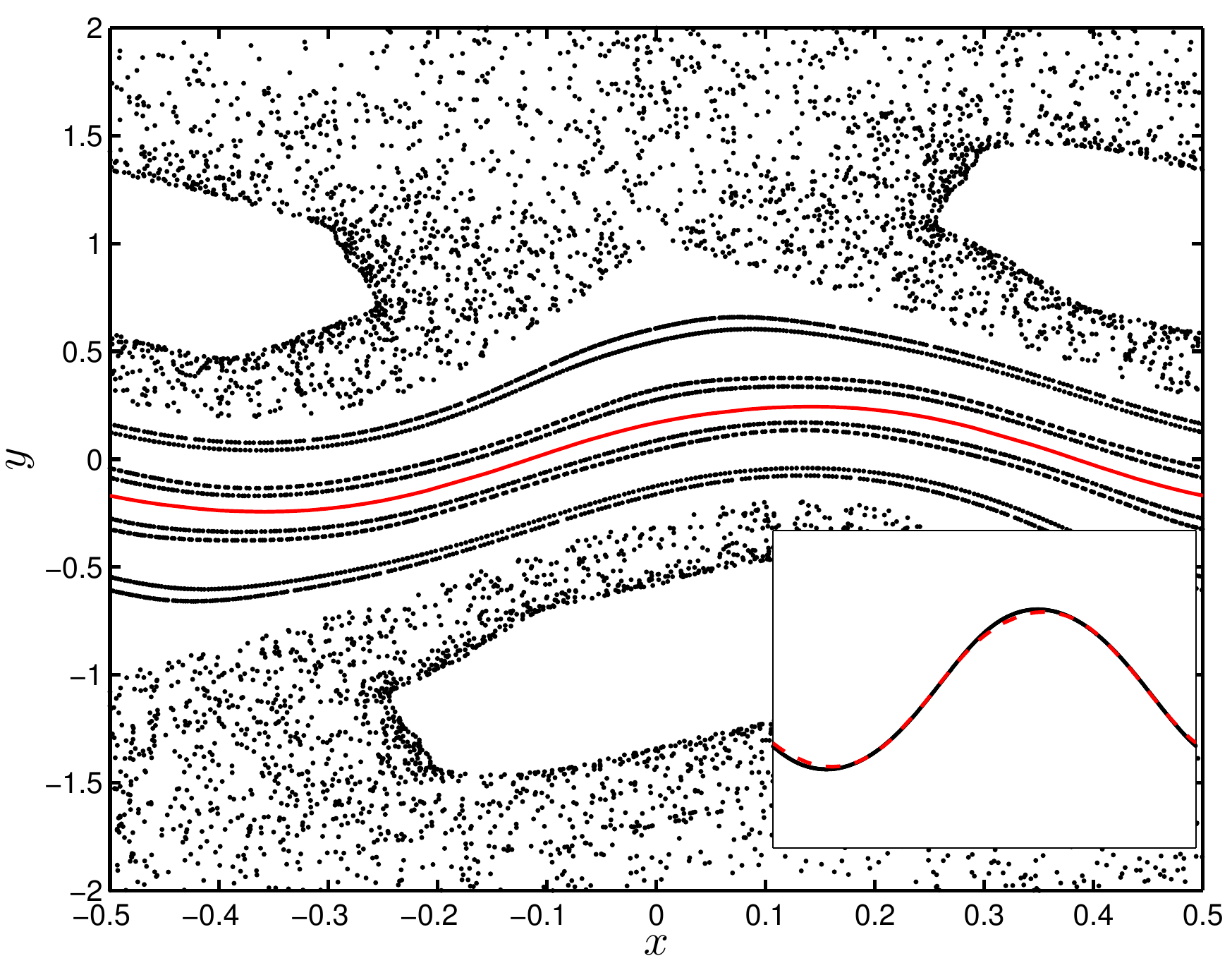}
\end{center}
\caption{The chaotic SNTM with parameters $a=0.27$, $b=0.38$.
The red curve shows the parabolic barrier computed from
$100$ iterations of SNTM. The inset compares the parabolic barriers
with the exact asymptotic barrier (black curve) obtained by $200$ iterations of the indicator points.}
\label{fig:sntm_chaotic} 
\end{figure}

The SNTM (\ref{eq:sntm}) becomes chaotic for parameters $a=0.27$,
$b=0.38$. The theory of indicator points still applies and gives
the exact asymptotic barrier for comparison. Figure \ref{fig:sntm_chaotic}
compares the computed parabolic barrier with the asymptotic shearless
barrier. The parabolic barrier is constructed from $100$ iterations of the SNTM while
the exact barrier is computed from $200$ iterations of the indicator point.

\subsection{Passive particles in mean-field coupled non-twist maps}

Following \cite{Diego_MF,Diego_MF_nt}, we consider the self-consistent mean field
interaction of $N$ coupled standard non-twist maps 
\begin{equation}
\begin{split} & x_{n+1}^{k}=x_{n}^{k}+a\left(1-\left(y_{n+1}^{k}\right)^{2}\right),\\
 & y_{n+1}^{k}=y_{n}^{k}-b_{n+1}\sin(2\pi x_{n}^{k}-\theta_{n}),
\end{split}
\label{eqn:mean_field_1}
\end{equation}
where $k=1,\ldots,N$ is an index for the particles and $n$ is the
iteration number. The variables $\theta_{n}$ and $b_{n}$ are given
by 
\begin{equation}
\begin{split} & \theta_{n+1}=\theta_{n}+\frac{1}{b_{n+1}}\frac{\partial\eta_{n}}{\partial\theta_{n}},\\
 & b_{n+1}=\sqrt{\left(b_{n}\right)^{2}+\left(\eta_{n}\right)^{2}}+\eta_{n},
\end{split}
\label{eqn:mean_field_2}
\end{equation}
where 
\begin{equation}
\eta_{n}=\sum_{i=1}^{n}\gamma_{i}\sin\left(x_{n}^{i}-\theta_{n}\right).\label{eqn:mean_field_3}
\end{equation}
We refer to the particles $x_{n}^{i}$ as \emph{active} particles
since they influence the mean field. The coefficients $\gamma_{i}$
are the coupling constants.  The mean field model \eqref{eqn:mean_field_1}-\eqref{eqn:mean_field_2} arises from studying vorticity defects in perturbations of parallel shear flow \cite{Diego_MF, Morrison_SWM}, and also has applications in one-dimensional beam plasmas \cite{Meiss_beam, Diego_MF}.

The full mean-field system is $2N$-dimensional, and we consider the behavior
of a \emph{passive} particle, whose non-autonomous evolution is given
by 
\begin{equation}
\begin{split} & x_{n+1}=x_{n}+a\left(1-y_{n+1}^{2}\right),\\
 & y_{n+1}=y_{n}-b_{n+1}\sin(2\pi x_{n}-\theta_{n}),
\end{split}
\label{eq:sntm_nonAut}
\end{equation}
where $b_{n}$ and $\theta_{n}$ are determined by the mean field
of active particles. The evolution of a passive particle is similar
to that of the SNTM considered in Section \ref{section:SNTM}, but
the parameters $b_{n}$ and $\theta_{n}$ change under each iteration
according to the mean field interaction of the active particles. When
the coupling constants $\gamma_{i}$ are zero, system (\ref{eq:sntm_nonAut})
coincides with the autonomous SNTM (\ref{eq:sntm}).

\begin{figure}[t!]
\begin{center}
\includegraphics[width=0.45\textwidth]{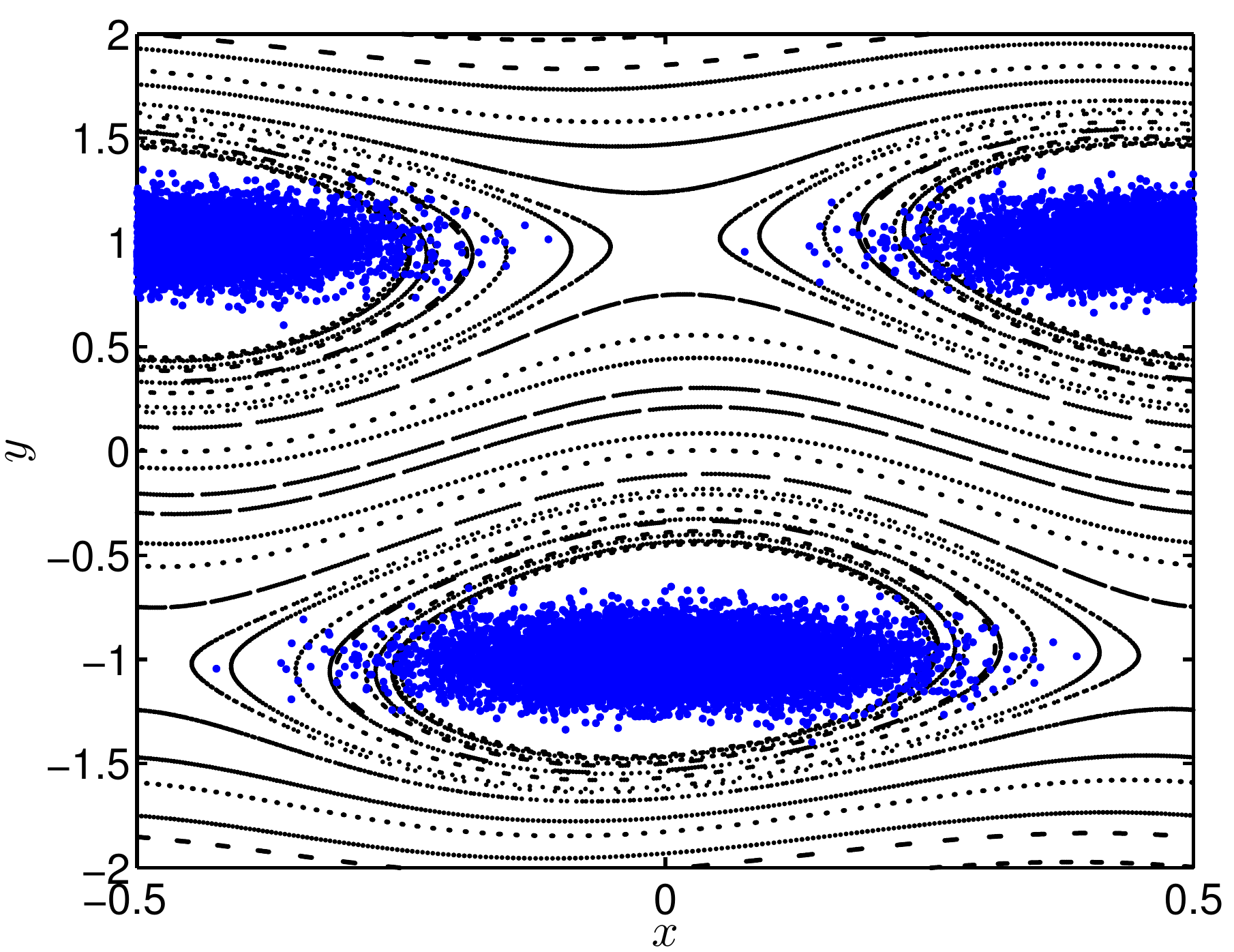} 
\includegraphics[width=0.45\textwidth]{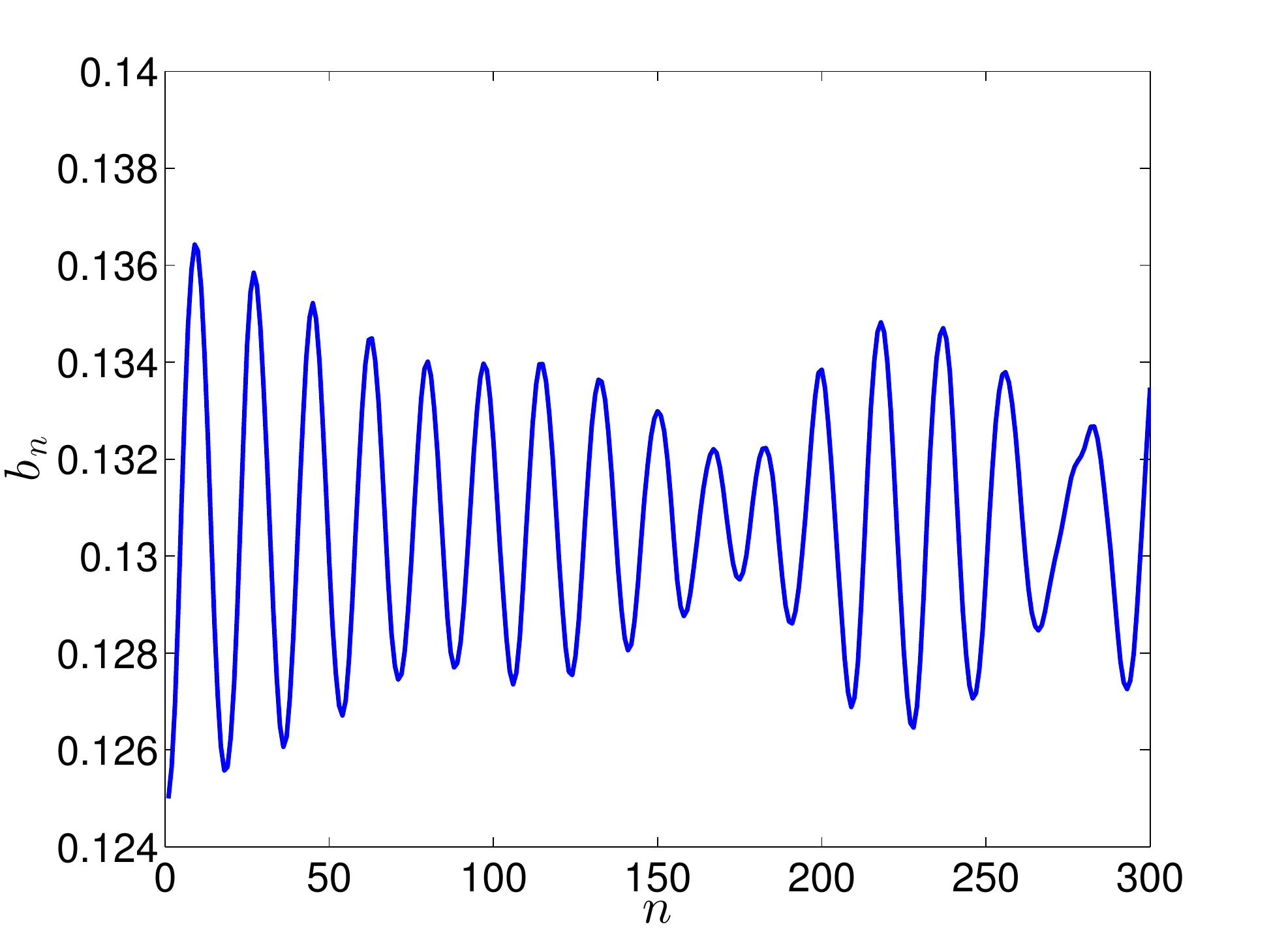}
\end{center}
\caption{Left: Initial conditions of the active particles. Right: Aperiodic
evolution of $b_{n}$.}
\label{fig:active_part_sntm} 
\end{figure}

We take $a=0.08$ and $b_{0}=0.125$ and $\theta_{0}=0.0$. The corresponding
dynamics for the SNTM (\ref{eq:sntm}) are integrable as described
in the previous section. With these initial parameters, we place $N=2\times10^{4}$
active particles localized near the islands (see Fig. \ref{fig:active_part_sntm})
and compute their mean field evolution. The coupling constants $\gamma_{i}$
are $2\times10^{-5}$ for all $i$. The evolution of the parameter
$b_{n}$ is shown in Figure \ref{fig:active_part_sntm}, and one thus
sees that the evolution of a passive particle is aperiodic with respect
to the iteration number.

\begin{figure}[H]
\begin{center}
\includegraphics[width=0.4\textwidth]{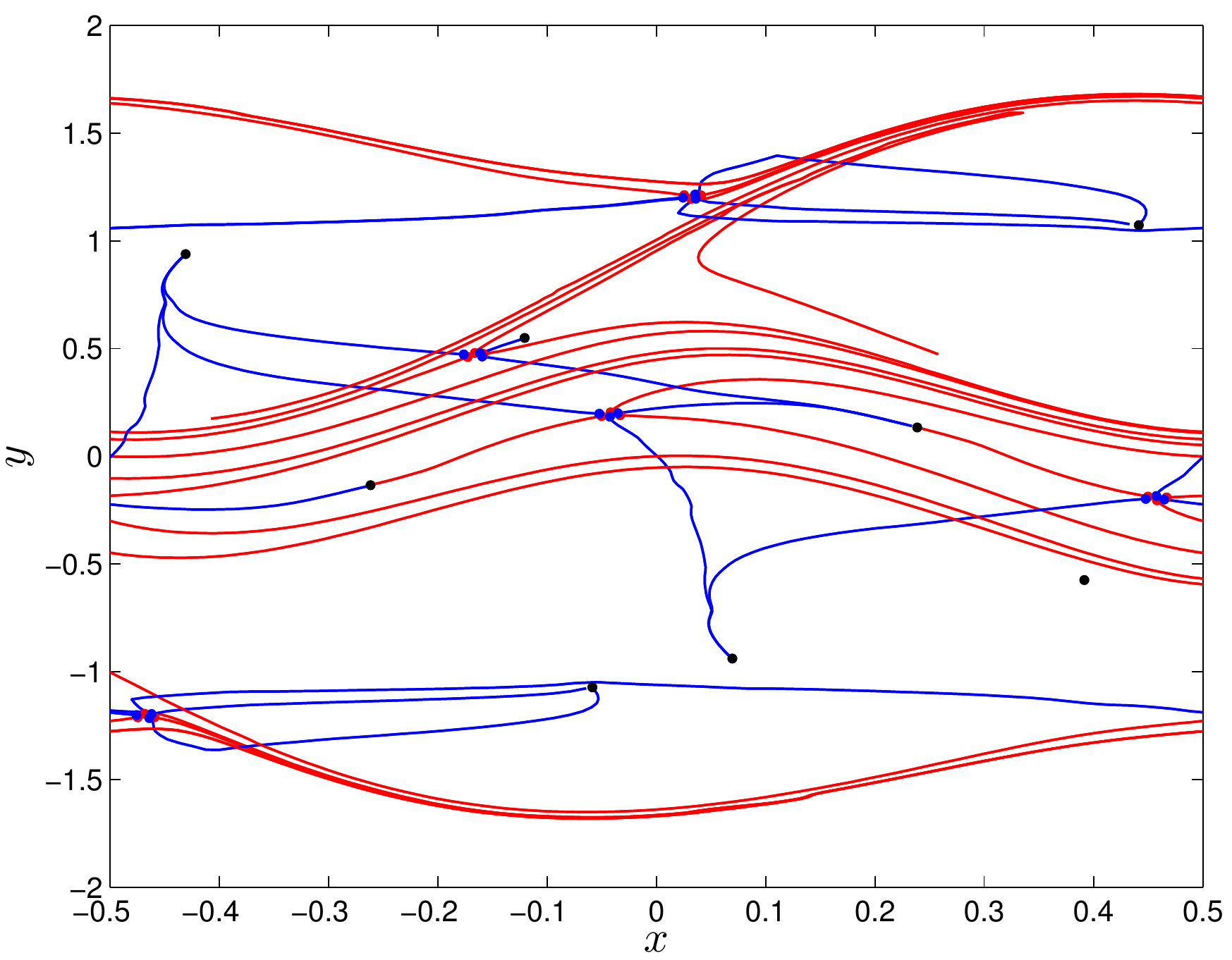}
\includegraphics[width=0.4\textwidth]{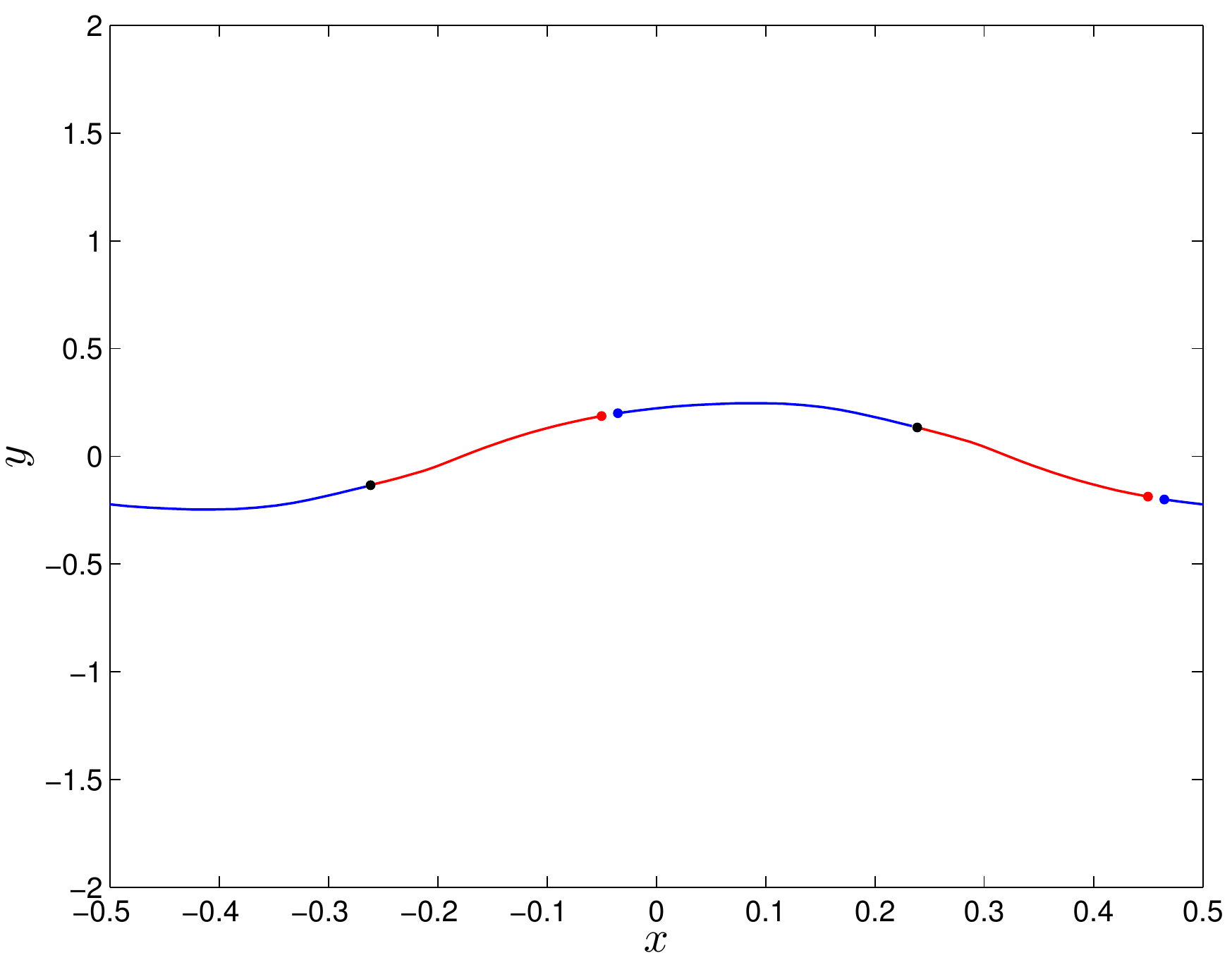} 
\end{center}
\caption{Left: Tensorlines for passive tracers in \eqref{eq:sntm_nonAut} strainlines (red)
and stretchlines (blue). The black dots mark the wedge singularities
where the tensorlines end. Right: Parabolic barrier as an alternating
sequence of tensorlines satisfying conditions P1-P2 of Definition
\ref{def:parabolic barriers}.}
\label{fig:tensorlines_sntm_naut} 
\end{figure}

With this setting, we compute all heteroclinic tensorlines using the
automated algorithm described in Section \ref{sec:methods}. Shown
in the left plot of Figure \ref{fig:tensorlines_sntm_naut}, the extracted
heteroclinic tensorline geometry is more complicated than what we
found for the SNTM. However, as seen in the right-side plot of the
figure, the final subset of connections satisfying conditions P1-P2
of Definition \ref{def:parabolic barriers} is similar to that of the
integrable system. This implies the persistence of a parabolic shearless barrier for a passive tracer
 in a self-consistent mean-field model.

\begin{figure}[H]
\begin{center}
\includegraphics[width=0.3\textwidth]{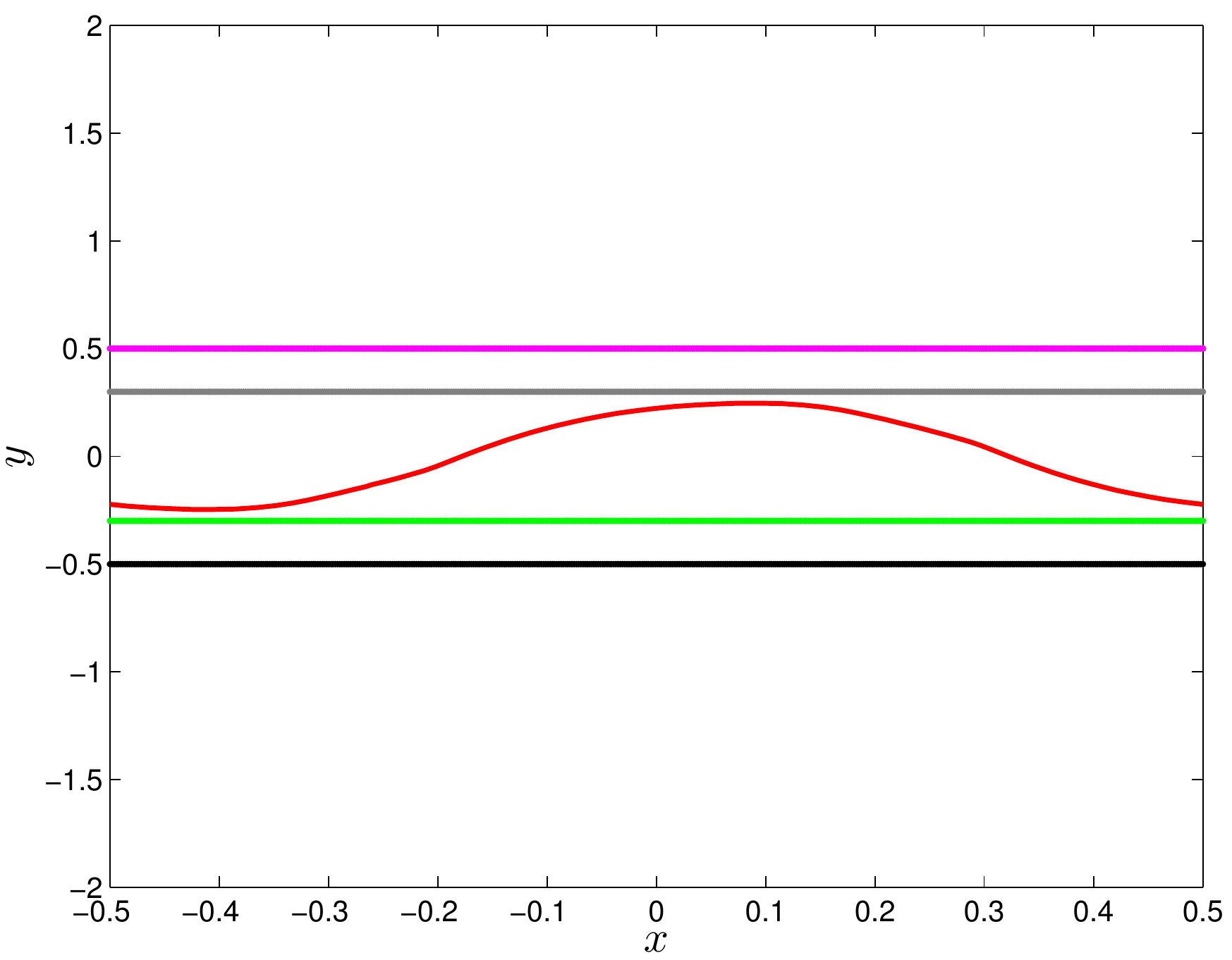}
\includegraphics[width=0.3\textwidth]{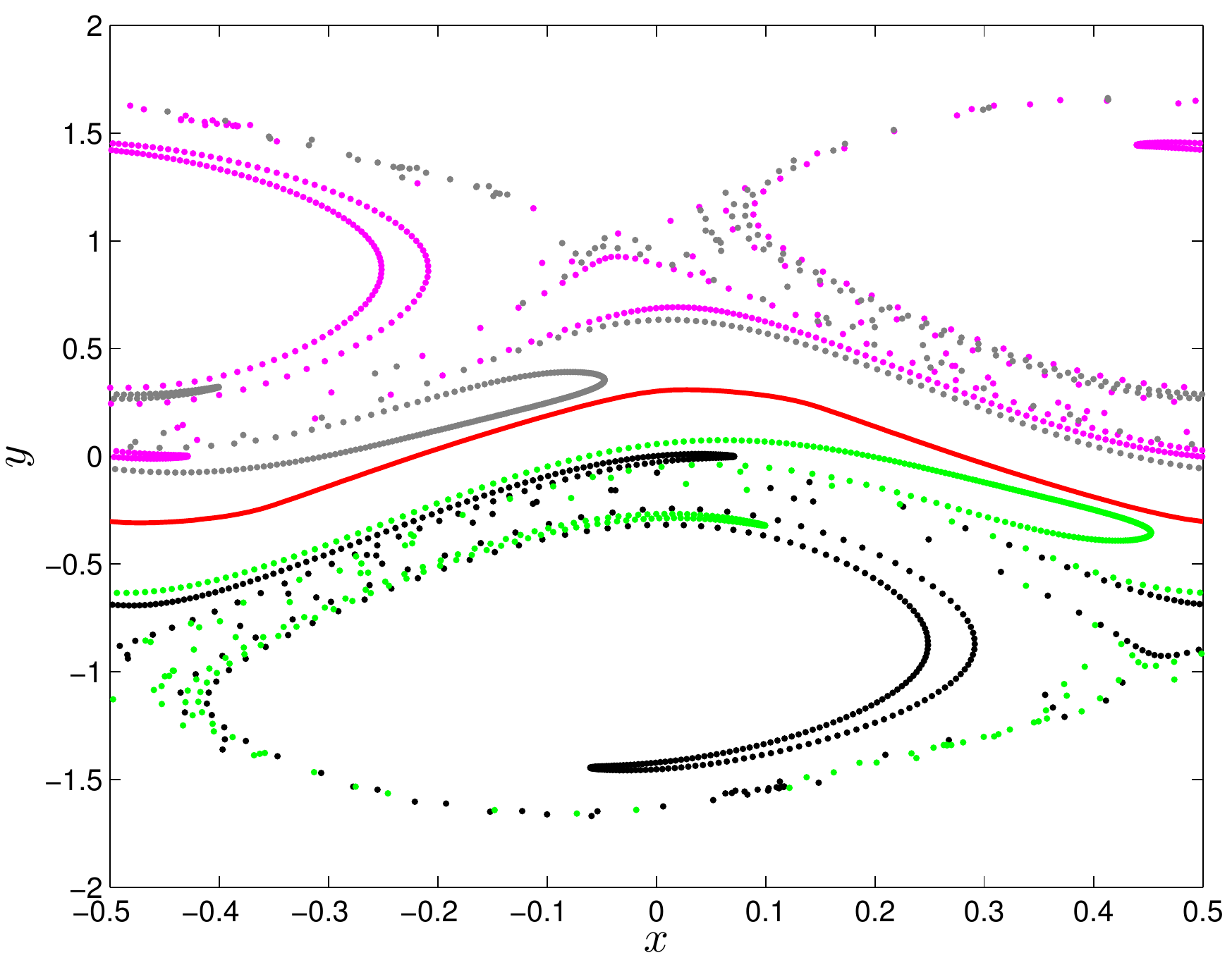}
\includegraphics[width=0.3\textwidth]{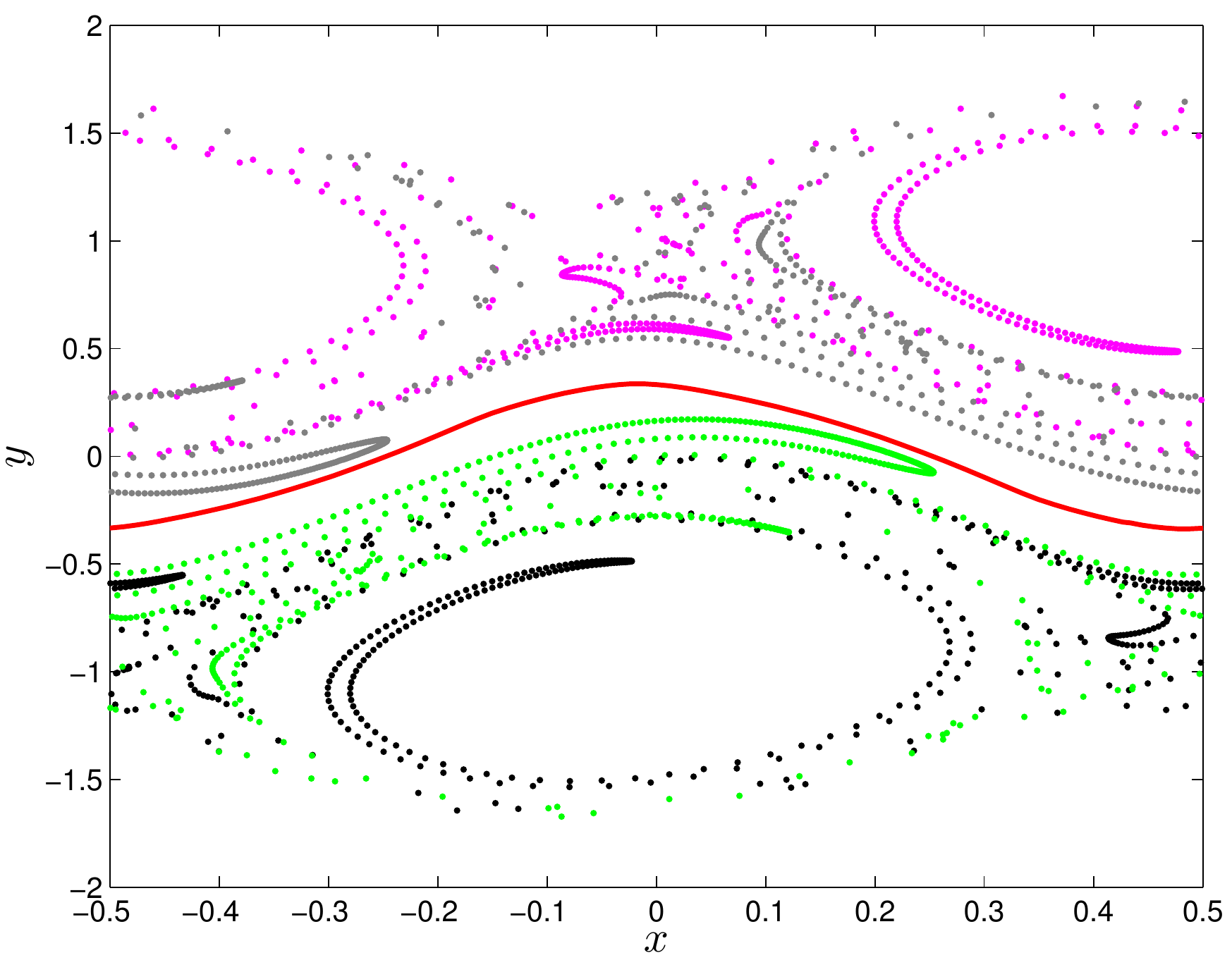}
\end{center}
\caption{Left: Parabolic barrier (red) and tracer particles (straight lines)
at the initial time. Advected images of the parabolic barrier and
tracer particles are shown after $50$ iterations (middle) and a $100$
iterations (right)}
\label{fig:partAdc_sntm_naut} 
\end{figure}

The evolution of tracers around the parabolic barrier is similar to
that shown in Figure \ref{fig:tracers_sntm_int}. Instead of presenting
the tracer evolution, however, we illustrate the role of the parabolic
barrier by placing two horizontal lines of particles above and two
below the parabolic barrier (cf. left plot of Fig. \ref{fig:partAdc_sntm_naut}).
The middle and right plots in the same figure show the advected images
of these lines after $50$ and $100$ iterations, respectively. We
conclude that despite the generally chaotic mixing prevalent in the
map, the extracted parabolic barrier provides a sharp and coherent
dividing surface that inhibits transport of passive particles.

\subsection{Bickley jet}

As our last example, we consider an idealized model of an eastward
zonal jet known as the Bickley jet \cite{Dfluid,Javier_bickley} in
geophysical fluid dynamics. This model consists of a steady background
flow subject to a time-dependent perturbation. The time-dependent
Hamiltonian for this model reads 
\begin{equation}
\psi(x,y,t)=\psi_{0}(y)+\psi_{1}(x,y,t),
\end{equation}
where 
\begin{equation}
\psi_{0}(y)=-UL\tanh\left(\frac{y}{L}\right),
\end{equation}
is the steady background flow and 
\begin{equation}
\psi_{1}(x,y,t)=UL\mathrm{sech}^{2}\left(\frac{y}{L}\right)\mbox{Re}\left[\sum_{n=1}^{3}f_{n}(t)\exp(\mbox{i}k_{n}x)\right],\label{eq:psi_1}
\end{equation}
is the perturbation. The constants $U$ and $L$ are characteristic
velocity and characteristic length scale, respectively. For the following
analysis, we apply the set of parameters used in \cite{Javier_bickley}:
\begin{equation}
U=62.66\;\mbox{ms}^{-1},\ \ L=1770\;\mbox{km},\ \ k_{n}=2n/r_{0},
\end{equation}
where $r_{0}=6371\;\mbox{km}$ is the mean radius of the earth.

\subsubsection{Quasiperiodic Bickley jet}
For $f_{n}(t)=\epsilon_{n}\exp(-\mbox{i}k_{n}c_{n}t)$, the time-dependent
part of the Hamiltonian consists of three Rossby waves with wave-numbers
$k_{n}$ traveling at speeds $c_{n}$. The amplitude of each Rossby
wave is determined by the parameters $\epsilon_{n}$. For small constant
values of parameters $\epsilon_{n}$, the Bickley jet is known to
have a closed, shearless jet core. In \cite{Javier_shear}, it is shown numerically
that this jet core is marked by a trench of the forward- and backward-time
FTLE fields. This finding is a consequence of temporal quasi-periodicity
of Rossby waves, which renders the the forward- and backward-time
dynamics as similar. In general, however, the time-dependence $f_{n}(t)$
can be any smooth signal \cite{geo_theory} with no particular recurrence.
We focus here on the existence of the shearless jet core under such
general forcing functions.

\begin{figure}[H]
\begin{center}
\includegraphics[width=0.6\textwidth]{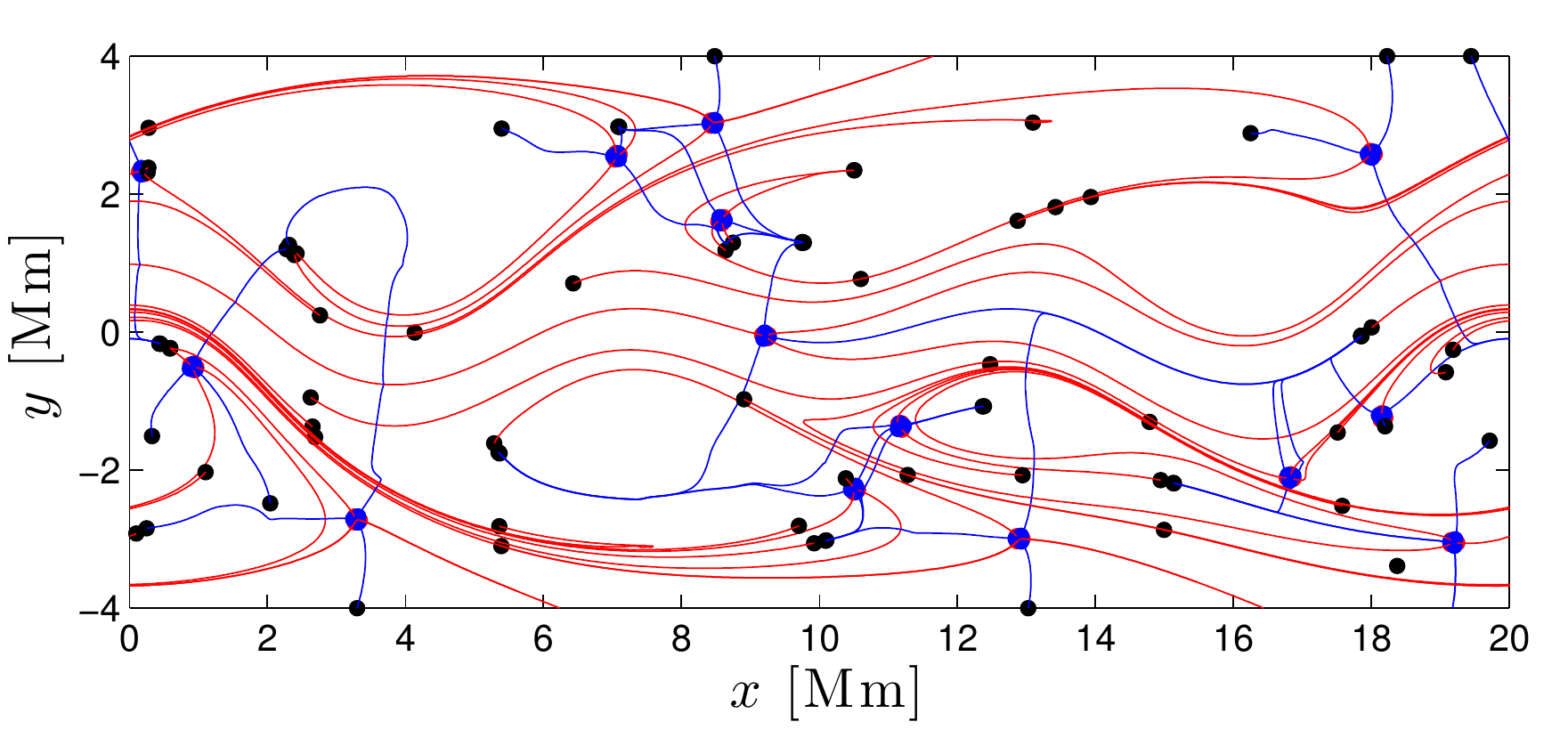}\\
\includegraphics[width=0.6\textwidth]{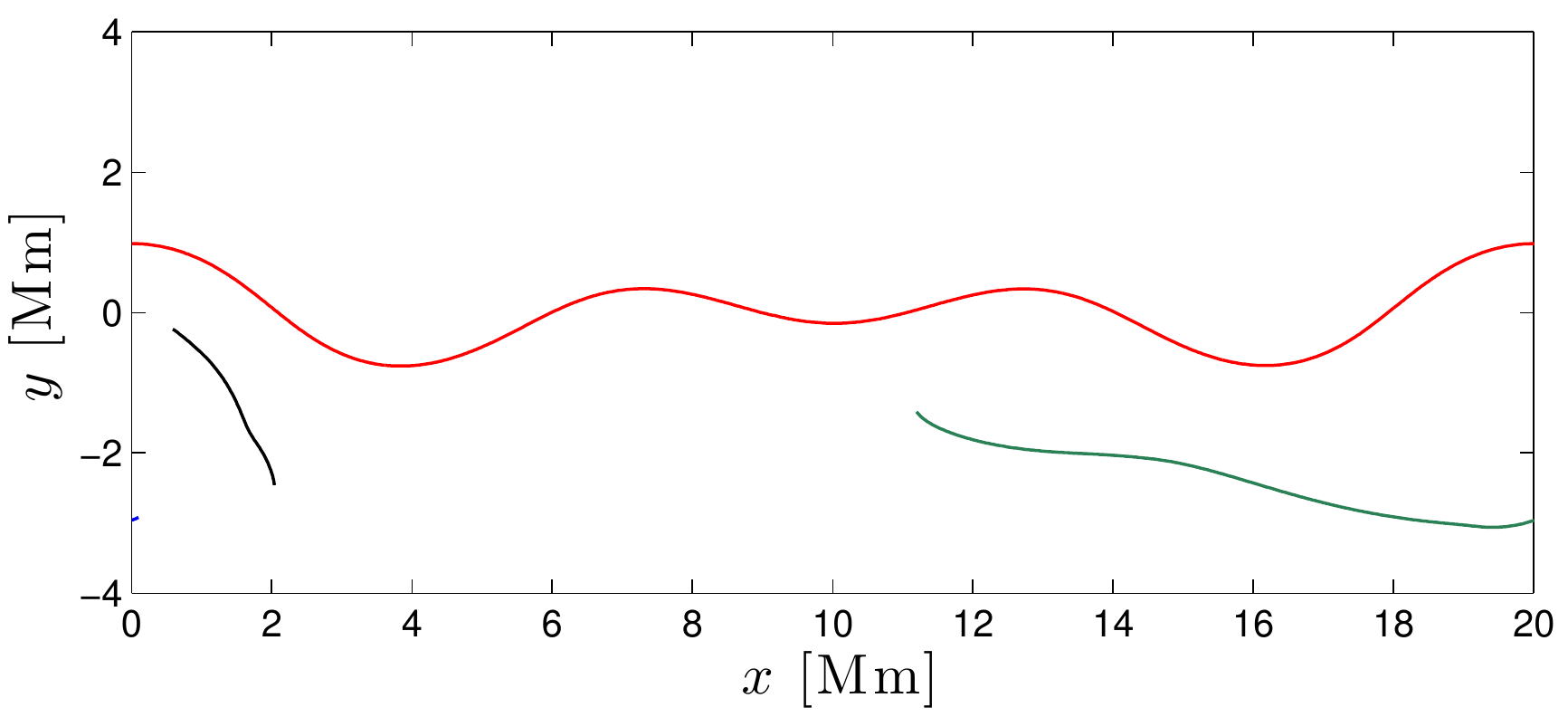}
\par\end{center}
\caption{Top: Tensorlines for the quasi-periodically forced Bickley jet: strainlines
(red) and stretchlines (blue). The black dots mark the wedge singularities
where the tensorlines end in while the blue dots mark the trisectors
where the tensorlines are initiated from. Bottom: Automatically extracted
parabolic barriers in the quasiperiodic Bickley jet.}
\label{fig:tensorlines_bickley_quasiPer} 
\end{figure}

First, however, we compare our results with those of \cite{Javier_shear}
for the quasi-periodic forcing $f_{n}(t)=\epsilon_{n}\exp(-\mbox{i}k_{n}c_{n}t)$,
with constant amplitudes $\epsilon_{1}=0.075$, $\epsilon_{2}=0.4$
and $\epsilon_{3}=0.3$. The top plot of Fig. \ref{fig:tensorlines_bickley_quasiPer}
shows automatically extracted heteroclinic tensorlines initiated from
trisectors and ending in wedges. Out of all these connections, three
satisfy conditions P1-P2 of Definition \ref{def:parabolic barriers}
and hence qualify as parabolic barriers (bottom plot of Fig. \ref{fig:tensorlines_bickley_quasiPer}). 

\begin{figure}[H]
\begin{center}
\includegraphics[width=0.31\textwidth]{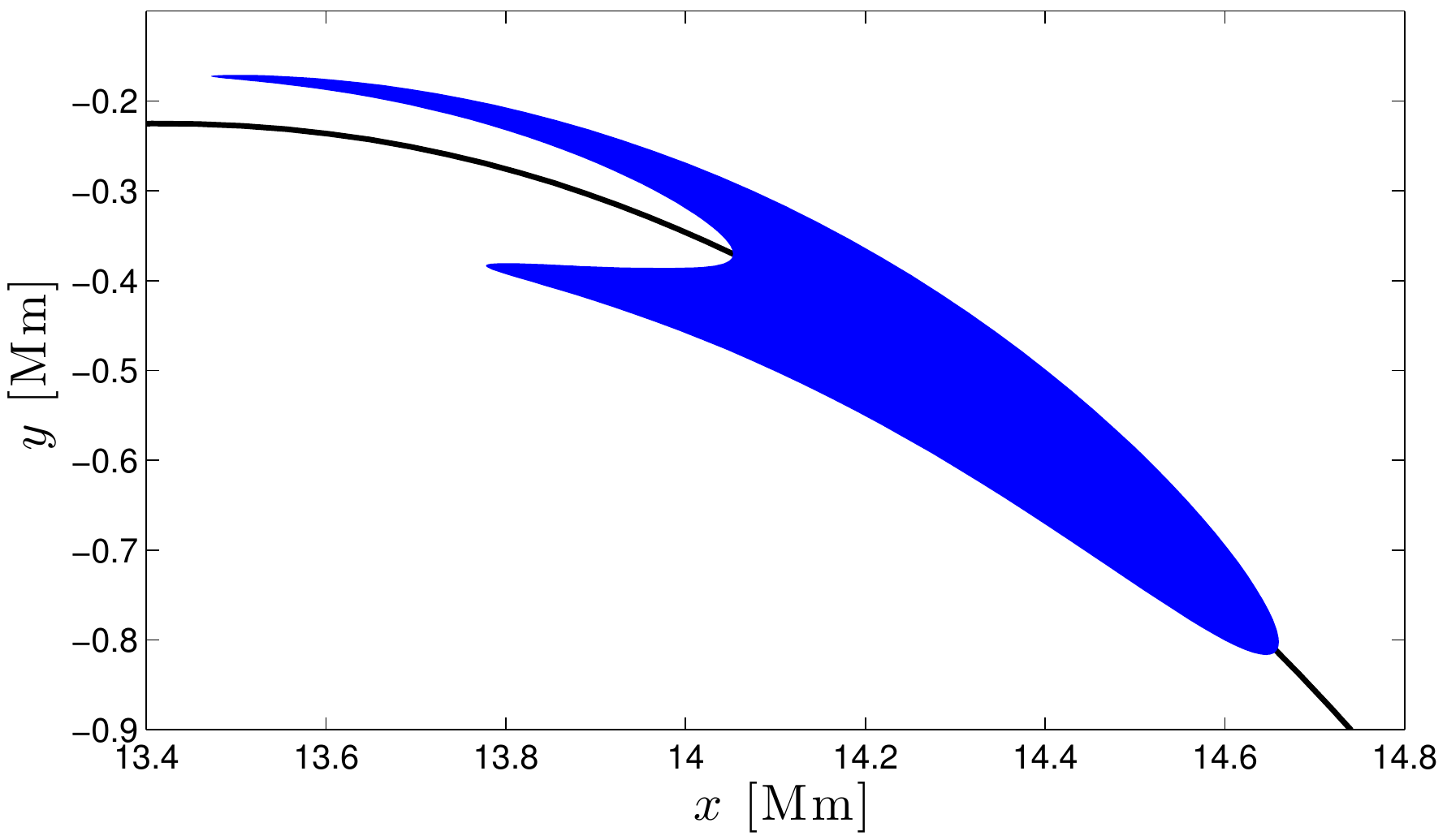}
\includegraphics[width=0.36\textwidth]{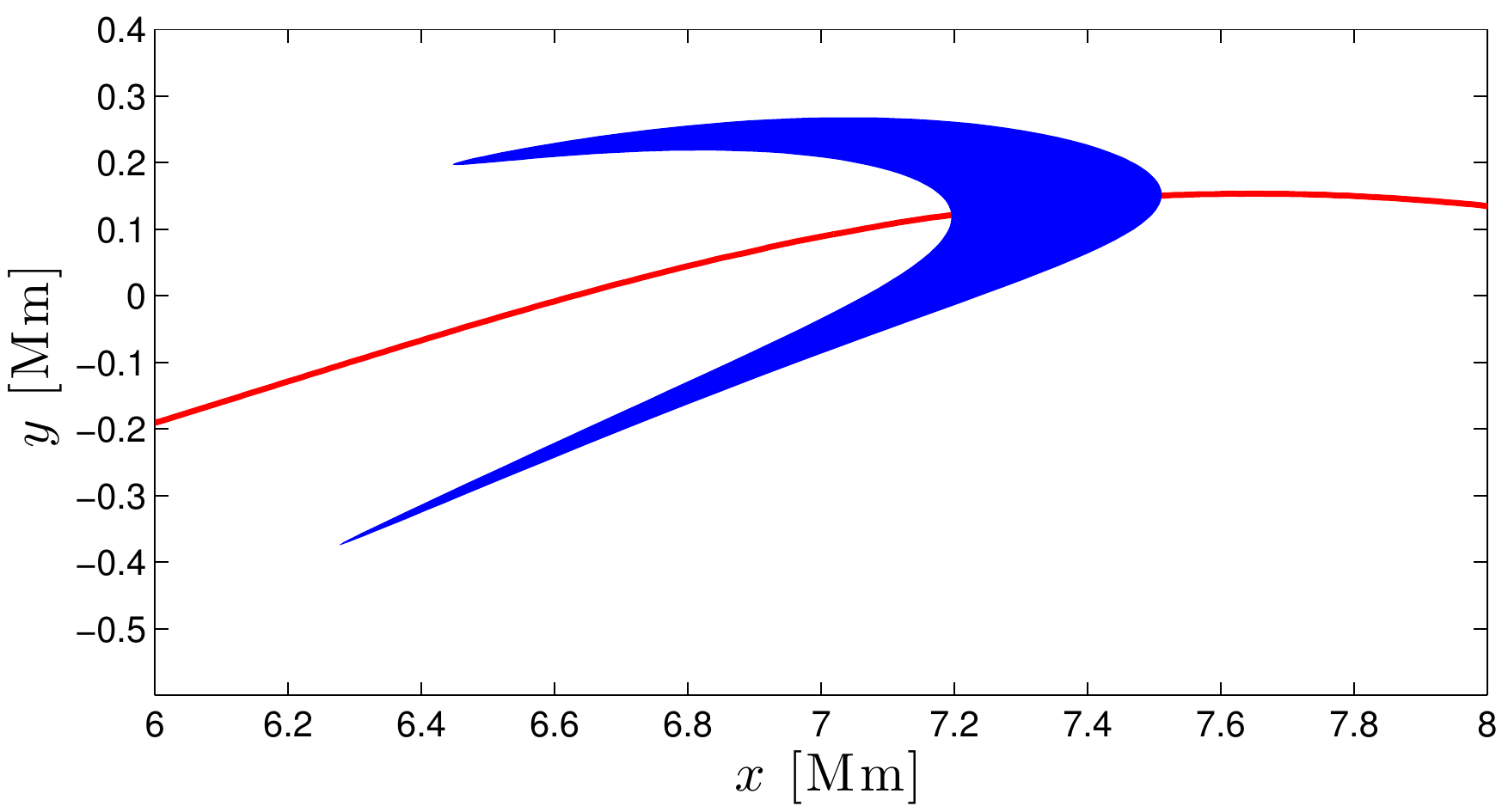}
\includegraphics[width=0.3\textwidth]{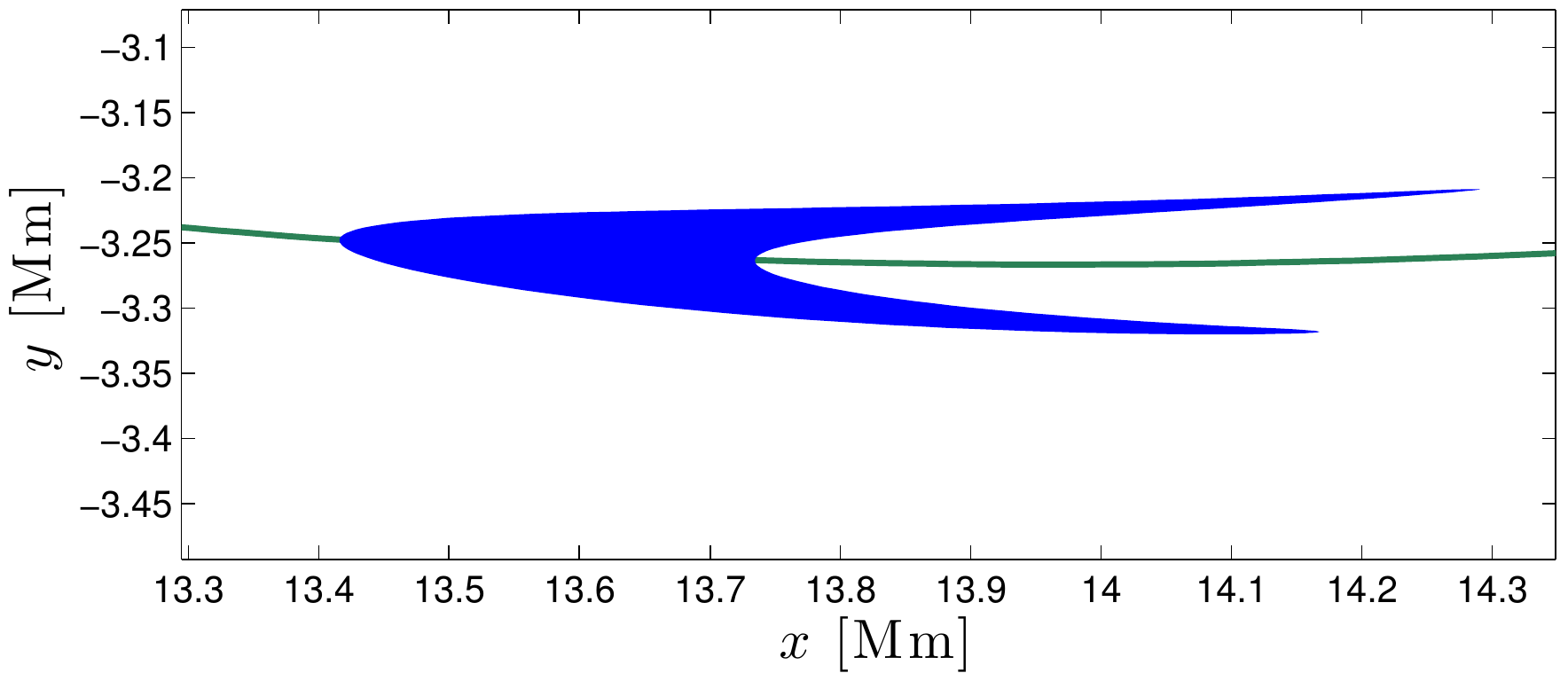}
\par\end{center}
\caption{The deformation of initially circular tracers (of radii $0.2$ Mm)
centered on the shearless curves after $11$ days. The color of the
curves correspond to those of Fig. \ref{fig:tensorlines_bickley_quasiPer}.}
\label{fig:tracers_bickley_qper} 
\end{figure}

The closed ($x$-periodic) parabolic barrier in red has also been
obtained in \cite{Javier_shear} as a trench of both the forward and
the backward FTLE field. The other two open parabolic barriers (blue
and black), however, have remained undetected in previous studies
to the best of our knowledge. Yet these open parabolic barriers do
serve as cores of smaller-scale jets, as demonstrated by the distinct
boomerang-shaped patterns developed by tracer blobs initialized along
them (see Fig. \ref{fig:tracers_bickley_qper}).

Such shearless material curves do not exist in the steady or
time-periodic counterpart of the Bickley jet, and thus perturbative
theories, such at KAM-type arguments, would not predict the existence
of such a jet core. Moreover, since these curves are not 
closed barriers separating the phase space they cannot be detected as almost-invariant coherent sets \cite{froyland10}.

\subsubsection{Chaotically forced Bickley jet}

\begin{figure}[H]
\begin{center}
\includegraphics[width=0.5\textwidth]{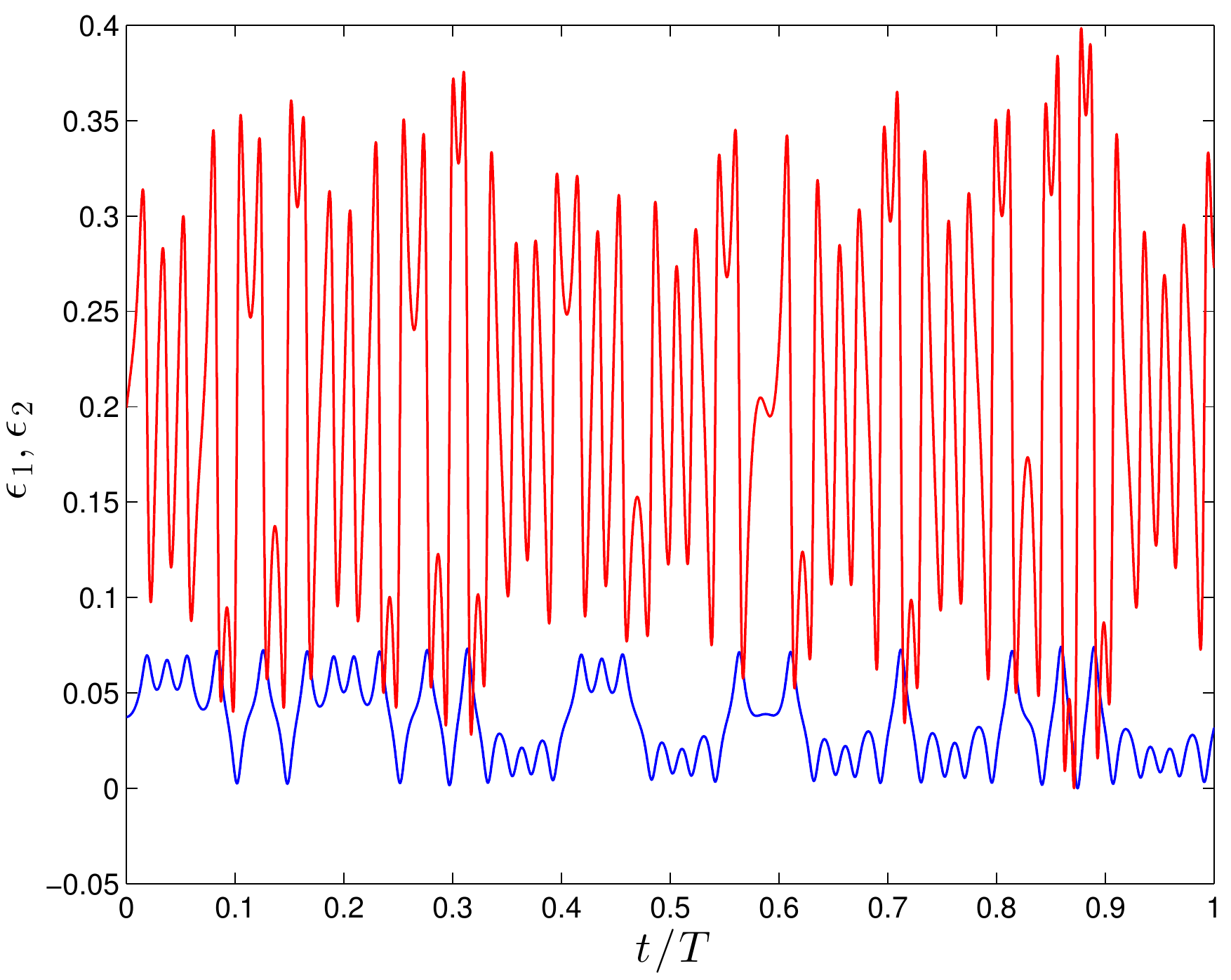} 
\end{center}
\caption{The chaotic signals $\epsilon_{1}$ (blue) and $\epsilon_{2}$ (red)
used as the amplitude of the forcing in equation (\ref{eq:psi_1}).
The integration time $T$ is 11 days.}
\label{fig:chaotic_signal} 
\end{figure}

To generate chaotic forcing for the Bickley jet, we let the forcing
amplitudes $\epsilon_{n}$ to be a chaotic signal for $n=1,2$. The
forcing amplitude $\epsilon_{3}=0.3$ remains constant. Figure \ref{fig:chaotic_signal},
shows the chaotic signals $\epsilon_{1}(t)$ and $\epsilon_{2}(t)$.

\begin{figure}[t!]
\begin{center}
\includegraphics[width=0.6\textwidth]{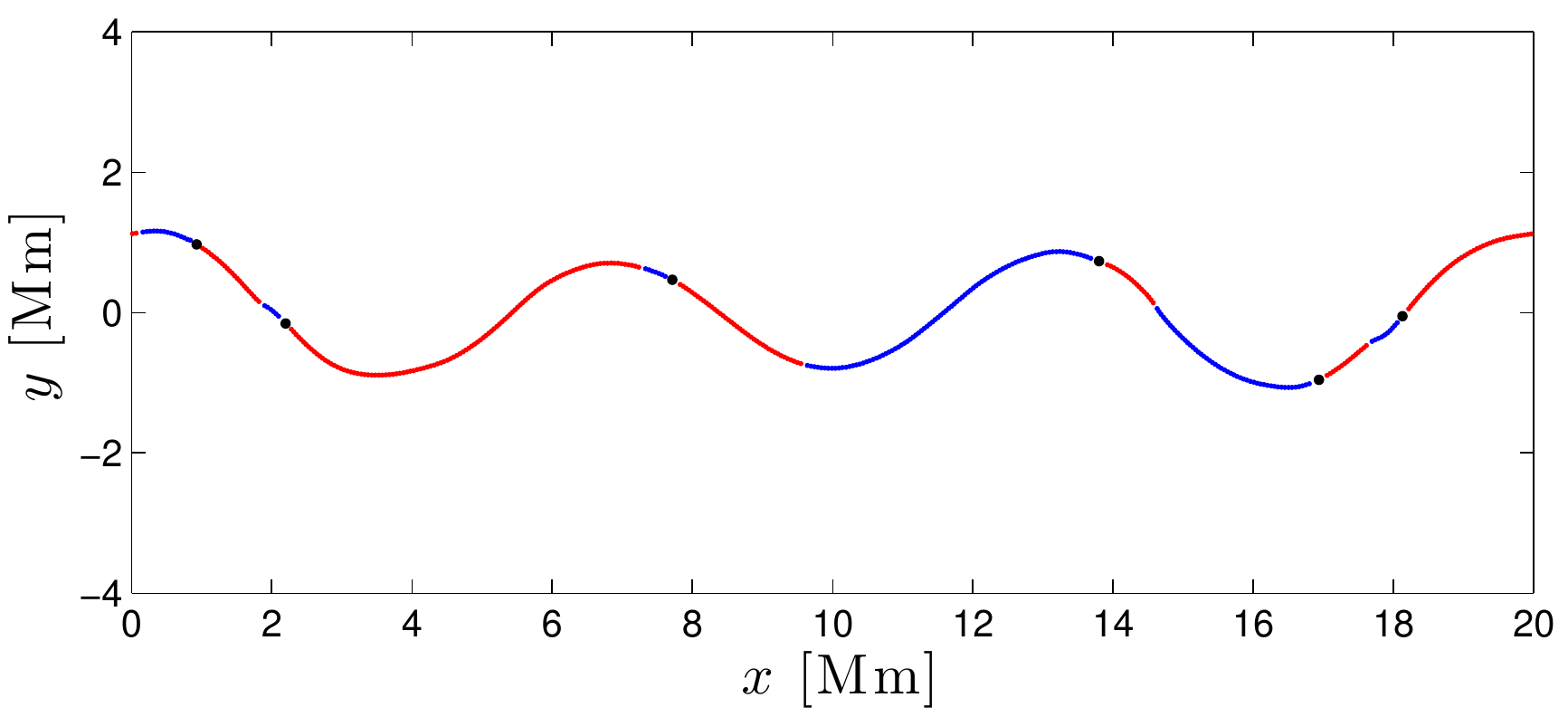}
\end{center}
\caption{The shearless curve for chaotically forced Bickley jet. The shearless
curve consists of alternating sequence of strainlines (red) and
stretchlines (blue). The wedge singularities are marked by black
dots.}
\label{fig:bickley_chaotic} 
\end{figure}

Figure \ref{fig:bickley_chaotic} shows the single parabolic barrier
obtained from the automated extraction procedure described in Section
\ref{sec:methods}. The additional open parabolic barriers found in 
the quasi-periodically forced case are, therefore, destroyed under
chaotic forcing.

\begin{figure}[t!]
\begin{center}
\includegraphics[width=0.45\textwidth]{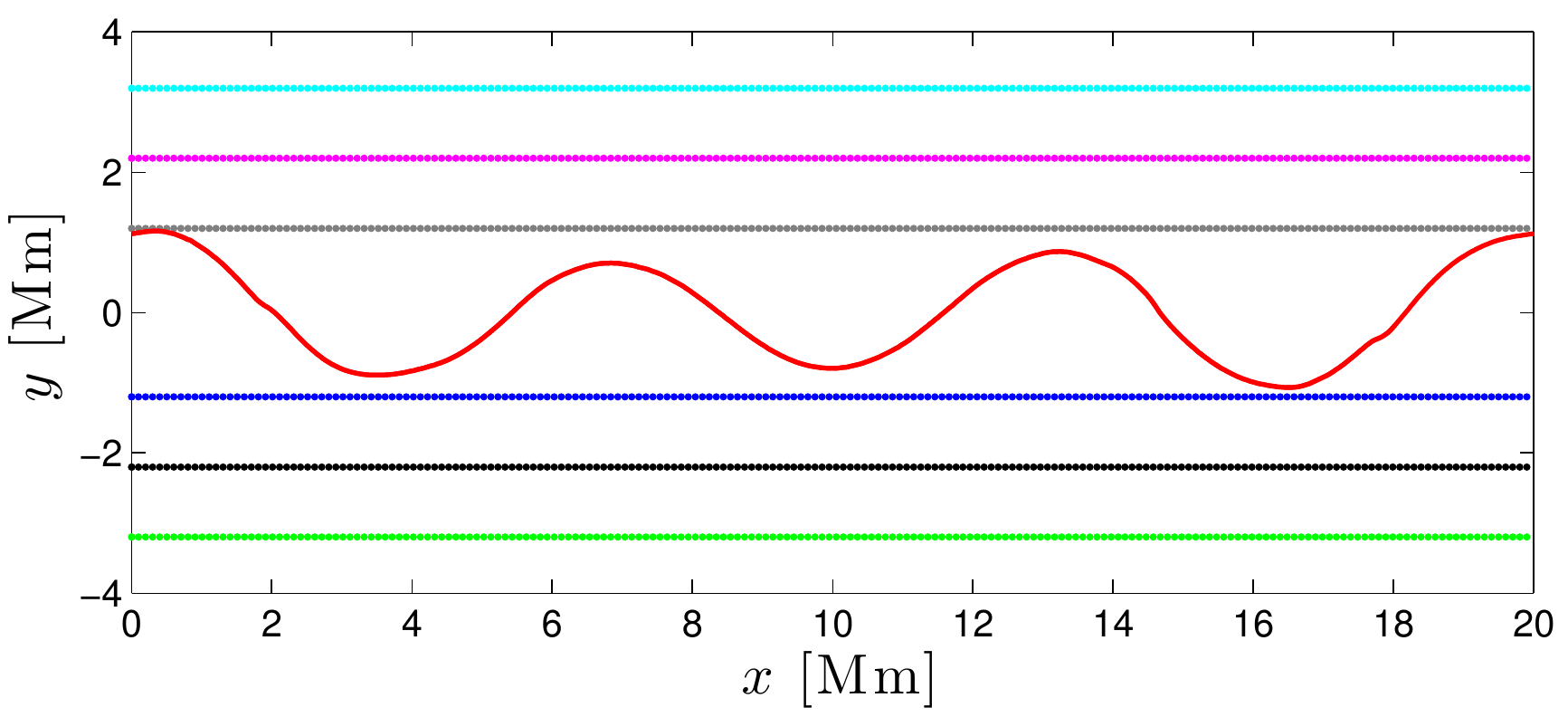}
\includegraphics[width=0.45\textwidth]{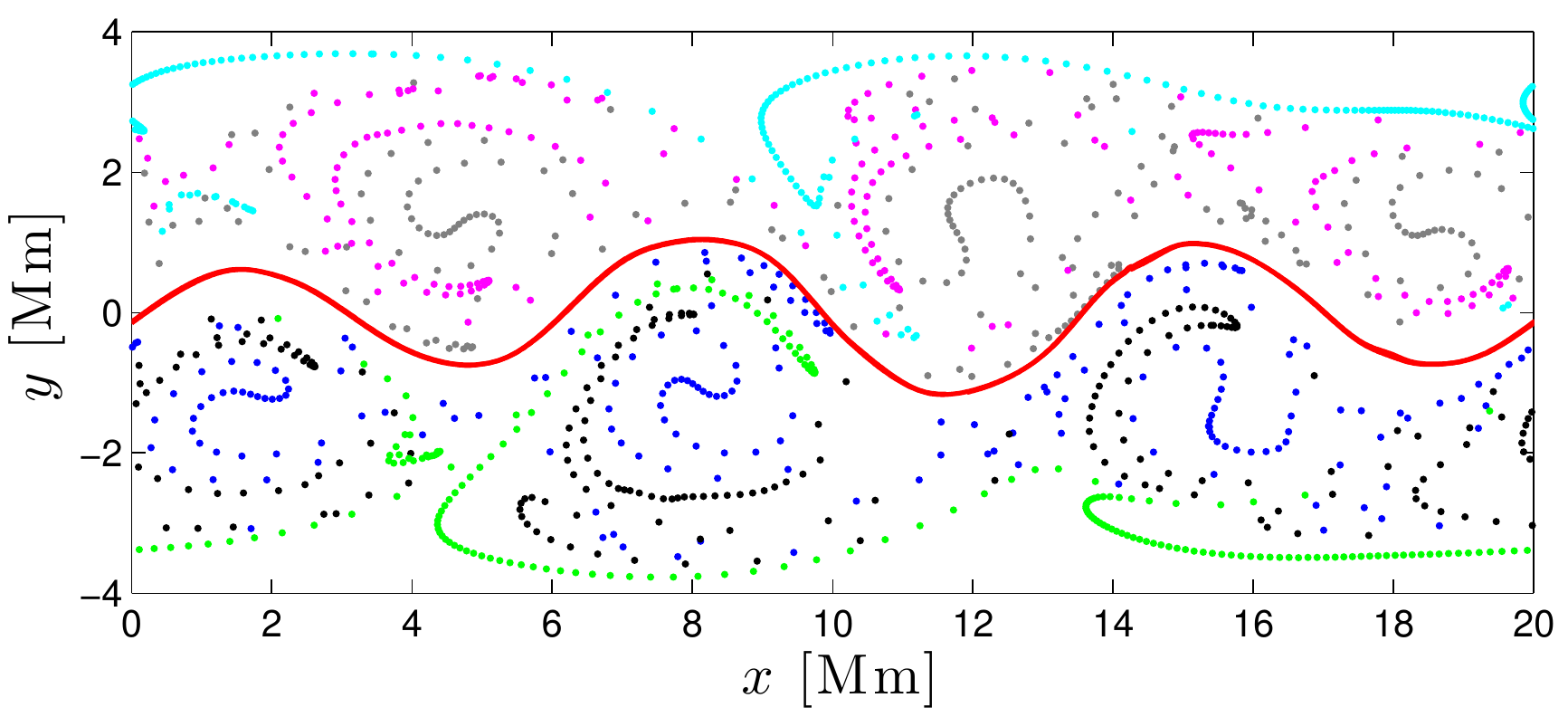}\\
\includegraphics[width=0.45\textwidth]{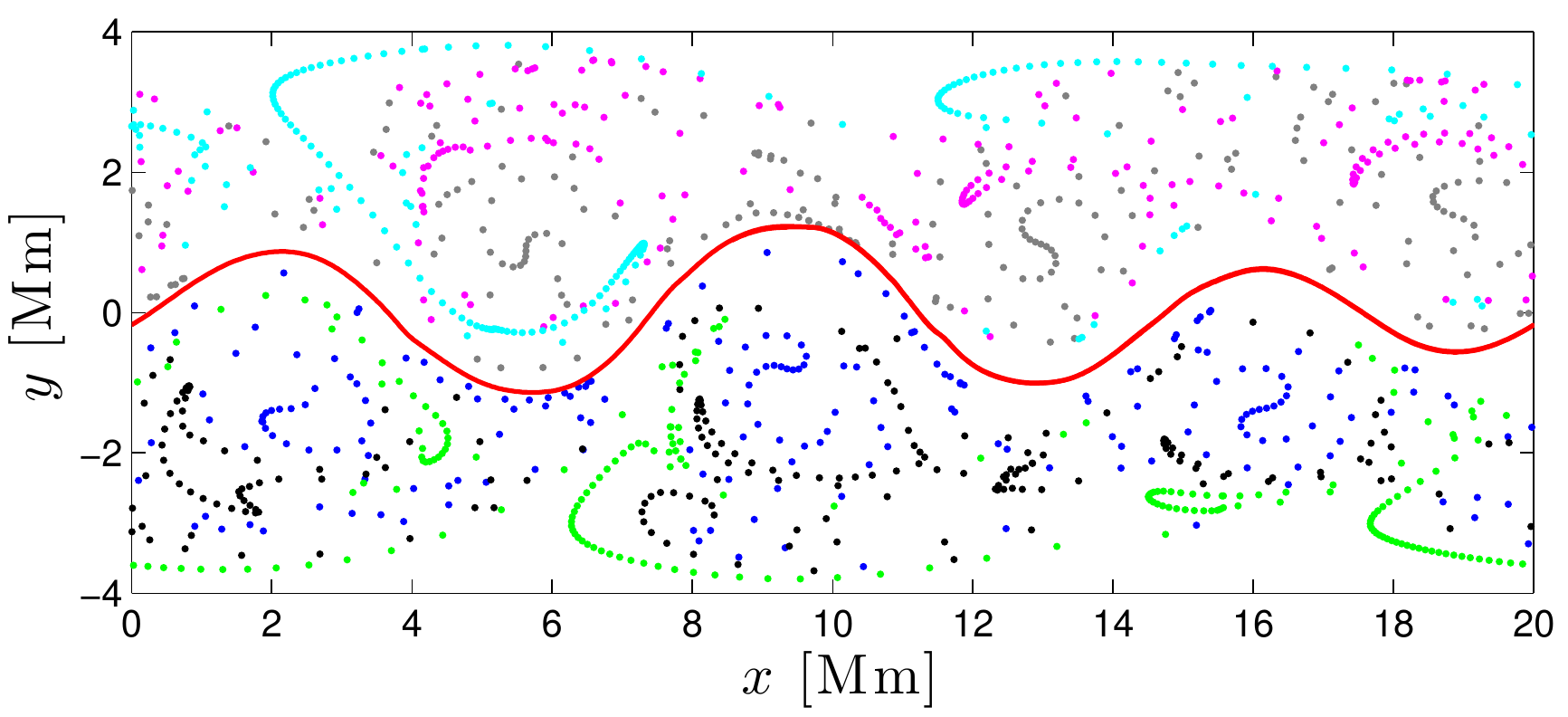}
\includegraphics[width=0.45\textwidth]{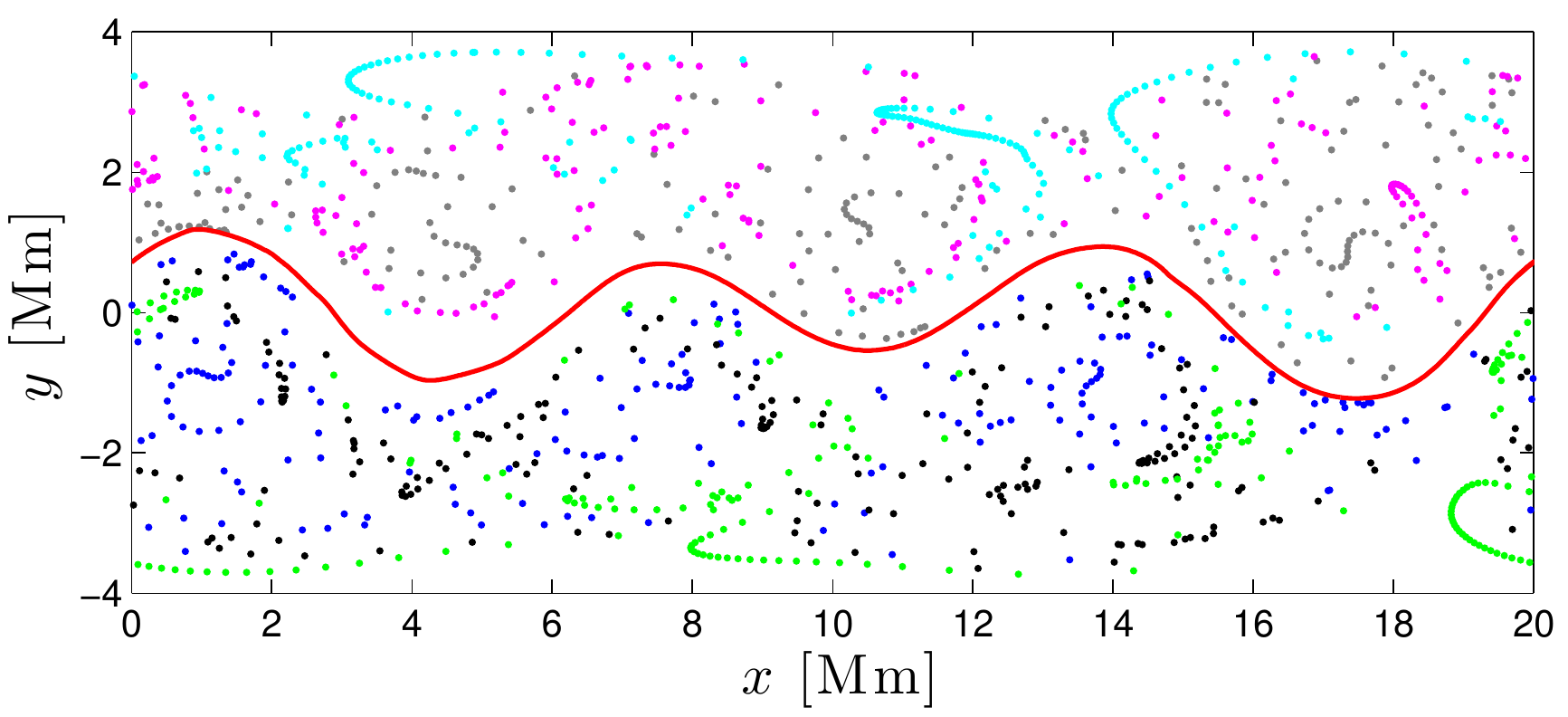}
\end{center}
\caption{Chaotically forced Bickley jet. The closed shearless curve (red) and
tracer particles (dots) at time $t=0$ (top left). Their advected
images are shown after 6 days (top right) 9 days (bottom left) and
11 days (bottom right).}
\label{fig:partAdv_bickley_chaotic} 
\end{figure}

The dynamic role of the remaining single barrier is illustrated in
Fig. \ref{fig:partAdv_bickley_chaotic}, where initially straight
lines of passive particles are advected for $6$, $9$ and $11$ days. 
Despite widespread chaotic mixing, the parabolic barrier
preserves its coherence, showing no stretching, folding, or smaller-scale
filamentation. Therefore, the extracted parabolic barrier is a sharp
separator between two invariant mixing regions. This shows that beyond
the almost-invariant sets located for the Bickley jet by set-theoretical
methods \cite{froyland09,froyland10}, actual invariant sets with sharp, coherent
boundaries also exist for the parameter values considered here.

\section{Conclusions}

We have developed a variational principle for shearless material lines
in two-dimensional, non-autonomous dynamical systems. Solutions to
this principle turn out to be composed of tensorlines of the Cauchy--Green
strain tensor. Locally most stretching or contracting tensorlines
staying away from singularities of the Cauchy--Green strain tensor
are found to be hyperbolic Lagrangian Coherent Structures (LCSs).
Thus, the present results give the first global variational description
of hyperbolic LCS as shearless material curves. 

By contrast, special chains of alternating tensorlines between Cauchy--Green
singularities define another class of shearless barriers, which we
call parabolic barriers (or parabolic LCSs). These barriers satisfy
variable-endpoint boundary conditions in the underlying Euler-Lagrange
equation, which make them exceptionally robust with respect to a broad
class of perturbations. This explains the broadly reported robustness
of shearless barriers observed in physical systems. 

We have devised an algorithm for the automated numerical detection
of parabolic barriers in two-dimensional unsteady flows. We illustrated
this algorithm on the standard non-twist map (SNTM), passive tracers
in mean-field coupled SNTMs and a model of the zonal jet (known as
the Bickley jet). For the SNTM, we showed that under increasing iterations,
our parabolic barrier converges to the exact shearless curve predicted
by the theory of indicator points. 

For the Bickley jet, we have recovered the results of \cite{Javier_shear}
on closed zonal jet cores under quasi-periodic forcing. We have also
found, however, other open jet cores in the same setting that were
not revealed by previous studies. A zonal jet was also detected in 
a chaotically forced Bickley jet.

While higher-dimensional shearless barriers have not yet been studied
extensively, the variational methods developed here should extend
to higher-dimensional flows. Such an extension of the concept of a
parabolic barrier appears to be possible via the approach developed
recently for elliptic and hyperbolic transport barriers in three-dimensional
unsteady flows \cite{LCS3D}.

\begin{appendices}

\section{Derivation of variable-endpoint boundary conditions for
the shearless variational principle}\label{bdrycond}

Note that
\begin{eqnarray}
\partial_{r^{\prime}}p & = & \frac{\left[2\left\langle r^{\prime},Cr^{\prime}\right\rangle \left\langle r^{\prime},r^{\prime}\right\rangle D-\left\langle r^{\prime},Dr^{\prime}\right\rangle \left\langle r^{\prime},r^{\prime}\right\rangle C-\left\langle r^{\prime},Dr^{\prime}\right\rangle \left\langle r^{\prime},Cr^{\prime}\right\rangle I\right]r^{\prime}}{\sqrt{\left\langle r^{\prime},Cr^{\prime}\right\rangle \left\langle r^{\prime},r^{\prime}\right\rangle }^{3}}\label{eq:partialp}\\
\nonumber 
\end{eqnarray}
Defining 
$$M:=\frac{2\left\langle r^{\prime},Cr^{\prime}\right\rangle \left\langle r^{\prime},r^{\prime}\right\rangle D-\left\langle r^{\prime},Dr^{\prime}\right\rangle \left\langle r^{\prime},r^{\prime}\right\rangle C-\left\langle r^{\prime},Dr^{\prime}\right\rangle \left\langle r^{\prime},Cr^{\prime}\right\rangle I}{\sqrt{\left\langle r^{\prime},Cr^{\prime}\right\rangle \left\langle r^{\prime},r^{\prime}\right\rangle }^{3}},$$
we have
\begin{equation}
\partial_{r^\prime}p=Mr^\prime.
\end{equation}
Any perturbation $h$ can be written as $h=h_\parallel+h_\perp$ where $h_\parallel$ and $h_\perp$ are, respectively, the tangential and orthogonal components of $h$ with respect to $r'$. Therefore, the boundary term in \eqref{eq:vari} can be written as
\begin{equation}
\langle \partial_{r'}p,h\rangle = \langle Mr',h_\perp\rangle. 
\end{equation}
Note that the term $\langle Mr',h_\parallel\rangle$ vanishes since $\langle Mr',r'\rangle=0$. 

Since $h_\perp$ is a scalar multiple of $\Omega r'$, the boundary term $\langle \partial_{r'}p,h\rangle$ vanishes if and only if $\langle Mr',\Omega r'\rangle=0$. Now expanding $r'$ in the Cauchy--Green eigenbasis as $r'=\alpha\xi_1+\beta\xi_2$, we get
\begin{equation}
\langle Mr',\Omega r'\rangle=\frac{(\alpha^2\lambda_1+\beta^2\lambda_2)(\alpha^2-\beta^2)(\lambda_2-\lambda_1)-\alpha^2\beta^2(\lambda_2-\lambda_1)^2}{(\alpha^2+\beta^2)^{1/2}(\alpha^2\lambda_1+\beta^2\lambda_2)^{3/2}},
\label{eq:normal_pert}
\end{equation}
where we used the fact that $C\xi_i=\lambda_i\xi_i$ for $i=1,2$. Without loss of generality, we may assume that the tangent vector $r'$ is normalized such that $\alpha^2+\beta^2=1$. 

Clearly if $\lambda_2=\lambda_1$, $\langle Mr',\Omega r'\rangle$ vanishes and so does the boundary term $\langle \partial_{r'}p,h\rangle$. By definition, the eigenvalues $\lambda_1$ and $\lambda_2$ only coincide at the Cauchy--Green singularities. For an incompressible flow, $C=I$ at the Cauchy--Green singularities since $\lambda_1\lambda_2=1$. This proves the condition \eqref{eq:variBC}.

Alternatively, assuming $\lambda_1\neq\lambda_2$, we find that $\langle Mr',\Omega r'\rangle=0$ if and only if
\begin{equation*}
\alpha=\pm\sqrt{\frac{\sqrt{\lambda_2}}{\sqrt{\lambda_1}+\sqrt{\lambda_2}}},\ \ \ \beta=\pm\sqrt{\frac{\sqrt{\lambda_1}}{\sqrt{\lambda_1}+\sqrt{\lambda_2}}}.
\end{equation*}
In other words, for the boundary term $\langle \partial_{r'}p,h\rangle$ to vanish, the tangent vectors $r'$ at the endpoints of $\gamma$ must satisfy
$$r^\prime=\sqrt{\frac{\sqrt{\lambda_2}}{\sqrt{\lambda_1}+\sqrt{\lambda_2}}}\;\xi_1\pm\sqrt{\frac{\sqrt{\lambda_1}}{\sqrt{\lambda_1}+\sqrt{\lambda_2}}}\;\xi_2.$$

The above linear combination of the Cauchy--Green eigenvectors is referred to as the \emph{shear vector field} \cite{geo_theory}. \emph{Shearlines}, i.e. the solution curves of the shear vector field, have been shown to mark boundaries of coherent regions of the phase space \cite{geo_theory, hadji2013, oceanicEddies_agulhas}, e.g., generalized KAM tori and coherent eddy boundaries.

Shear vector fields, however, do not result in shearless transport barriers; in fact, they are local maximizers of Lagrangian shear \cite{geo_theory}. 

\section{Equivalent formulation of the shearless variational principle}\label{app:B}

\noindent With the shorthand notation 
\begin{equation}
A(r,r')=\langle r',C(r)r'\rangle,\quad B(r')=\langle r',r'\rangle,\quad G(r,r')=\langle r',D(r)r'\rangle,\label{eq:AB}
\end{equation}
$P$ can be rewritten as 
\begin{equation}
P(\gamma)=\frac{1}{\sigma}\int_{0}^{\sigma}p(r,r^{\prime})\,\mathrm{d}s=\frac{1}{\sigma}\int_{0}^{\sigma}\frac{G(r,r')}{\sqrt{A(r,r')B(r')}}\,\mathrm{d}s,
\end{equation}
and its Euler--Lagrange equations \eqref{eq:EL} can be re-written
as 
\begin{equation}
\partial_{r}\frac{G}{\sqrt{AB}}-\frac{d}{\mathrm{d}s}\partial_{r'}\frac{G}{\sqrt{AB}}=0.\label{eq:oneway}
\end{equation}
Note that 
\begin{equation}
\partial_{r}\frac{G}{\sqrt{AB}}=\frac{\partial_{r}G}{\sqrt{AB}}-\frac{G\left(B\partial_{r}A+A\partial_{r}B\right)}{2\sqrt{AB}^{3}},\quad\partial_{r^{\prime}}\frac{G}{\sqrt{AB}}=\frac{\partial_{r^{\prime}}G}{\sqrt{AB}}-\frac{G\left(B\partial_{r^{\prime}}A+A\partial_{r^{\prime}}B\right)}{2\sqrt{AB}^{3}},.\label{eq:ident}
\end{equation}

Since the integrand of $P(\gamma)$ has no explicit dependence on
the parameter $s$, Noether's theorem \cite{Gelfand-Fomin-00} guarantees
the existence of a first integral for \eqref{eq:oneway}. This integral
can be computed as 
\begin{equation}
I=\frac{G}{\sqrt{AB}}-\left\langle r',\partial_{r'}\frac{G}{\sqrt{AB}}\right\rangle =\frac{G}{\sqrt{AB}}=I_{0}=\mathrm{const},\label{eq:firstint}
\end{equation}
where we have used the specific form of the functions $A$ and $B$
from \eqref{eq:AB}, as well as the second equation from \eqref{eq:ident}.

With the notation $\mu=I_{0}$ , we therefore have the identity
\begin{equation}
G(r(s),r'(s))\equiv\mu\sqrt{A(r(s),r'(s))B(r'(s))}\label{eq:A=00003DDB}
\end{equation}
on any solution \eqref{eq:oneway} for some appropriate value of the
positive constant $\mu>0$ . 

We use the identity \eqref{eq:A=00003DDB} to rewrite the expressions
\eqref{eq:ident} as 
\begin{equation}
\partial_{r}\frac{G}{\sqrt{AB}}=\frac{1}{\sqrt{AB}}\partial_{r}\left[G-\mu\sqrt{AB}\right],\quad\partial_{r'}\frac{G}{\sqrt{AB}}=\frac{1}{\sqrt{AB}}\partial_{r^{\prime}}\left[G-\mu\sqrt{AB}\right].\label{eq:ident-1}
\end{equation}
We also introduce a rescaling of the independent variable $s$ in
equation \eqref{eq:oneway} via the formula 
\begin{equation}
\frac{\mathrm{d}\tau}{\mathrm{d}s}=\sqrt{A(r(s),r'(s))B(r'(s))},\label{eq:resc1}
\end{equation}
which, by the chain rule, implies 
\begin{equation}
\sqrt{A(r(s),r'(s))B(r'(s))}=\frac{1}{\sqrt{A(r(\tau),\dot{r}(\tau))B(\dot{r}(\tau))}},
\label{eq:resc2}
\end{equation}
with the dot referring to differentiation with respect to the new
variable $\tau$. Note that $\sqrt{A(r(s),r'(s))B(r'(s))}$ is non-vanishing
on smooth curves with well-defined tangent vectors, and hence the
change of variables \eqref{eq:resc1} is well-defined.

After the $s\mapsto\tau$ rescaling and the application of \eqref{eq:resc2},
the expressions in \eqref{eq:ident-1} imply 
\begin{eqnarray}
\partial_{r}\frac{G(r,r^{\prime})}{\sqrt{A(r,r^{\prime})B(r^{\prime})}} & = & \frac{\partial_{r}\left[G(r,\dot{r})-\mu\sqrt{A(r,\dot{r})B(\dot{r})}\right]}{\sqrt{A(r(\tau),\dot{r}(\tau))B(\dot{r}(\tau))}}\nonumber \\
\\
\frac{\mathrm{d}}{\mathrm{d}s}\partial_{r'}\frac{G(r,r^{\prime})}{\sqrt{A(r,r^{\prime})B(r^{\prime})}} & = & \frac{\frac{\mathrm{d}}{\mathrm{d}\tau}\partial_{\dot{r}}\left[G(r,\dot{r})-\mu\sqrt{A(r,\dot{r})B(\dot{r})}\right]}{\sqrt{A(r(\tau),\dot{r}(\tau))B(\dot{r}(\tau))}}.\nonumber \\
\end{eqnarray}
 Based on these identities, equation \eqref{eq:oneway} can be re-written
as 
\begin{equation}
\frac{1}{\sqrt{A(r(\tau),\dot{r}(\tau))B(\dot{r}(\tau))}}\left\{ \partial_{r}\left[G(r,\dot{r})-\mu\sqrt{A(r,\dot{r})B(\dot{r})}\right]-\frac{\mathrm{d}}{\mathrm{d}\tau}\partial_{\dot{r}}\left[G(r,\dot{r})-\mu\sqrt{A(r,\dot{r})B(\dot{r})}\right]\right\} =0.\label{eq:oneway-3}
\end{equation}

Since $1/\sqrt{A(r(\tau),\dot{r}(\tau))B(\dot{r}(\tau))}$ is non-vanishing
we obtain from \eqref{eq:oneway-3} that all solutions of \eqref{eq:oneway}
must satisfy the Euler--Lagrange equation derived from the Lagrangian
\begin{equation}
\mathcal{H}_{\mu}(r,\dot{r})=\frac{1}{2}\left[G(r,\dot{r})-\mu\sqrt{A(r,\dot{r})B(\dot{r})}\right].\label{eq:energy-1}
\end{equation}
Therefore, all stationary functions of the functional $P$ are also
stationary functions of the function for an appropriate value of $\mu$.
This value of $\mu$ can be determined from formula \eqref{eq:A=00003DDB},
which also shows that the corresponding stationary functions of $\mathcal{H}_{\mu}$
all satisfy 
\begin{equation}
\langle\dot{r}(\tau),D(r(\tau))\dot{r}(\tau)\rangle=\mu\sqrt{A(r,\dot{r})B(\dot{r})}.\label{eq:const2}
\end{equation}
For $\mu=0$, these solutions are null-geodesics of the Lorentzian
metric \eqref{eq:Dmetric} induced by the tensor $D$.

Conversely, assume that $r(\tau)$ satisfies both equations \eqref{eq:oneway-3}
and \eqref{eq:const2}. Reversing the steps leading to \eqref{eq:const2},
and employing the inverse rescaling of the independent variable as,
\begin{equation}
\frac{\mathrm{d}s}{\mathrm{d}\tau}=\sqrt{A(r(\tau),\dot{r}(\tau))B(\dot{r}(\tau))},
\end{equation}
we obtain that any rescaled solution $r(s)$ is also a solution of
the Euler--Lagrange equation \eqref{eq:oneway}. Therefore, each solution
of \eqref{eq:oneway-3} lying in the zero energy surface $\mathcal{H}_{\mu}(r,\dot{r})=0$
is also a stationary curve of the functional $P(\gamma)$, lying on
the energy surface $I(r,r')=\mu$, and hence satisfying the identity
\eqref{eq:A=00003DDB}.

\section{Tensorline singularities, heteroclinic tensorlines, and
their numerical detection }\label{app:C}

Here we briefly review some relevant aspects of tensorline geometry near
singularities of a symmetric tensor field \cite{Ten_lines_2D,Tricoche}.

\subsection{Tensorline singularities}

Singularities of tensorlines, such as the tensorlines of the Cauchy--Green
strain tensor, are points where the tensor field becomes the identity
tensor, and hence ceases to admit a well-defined pair of eigenvectors.
As a consequence, tensorlines, as curves tangent to $\xi_{1}$ and
$\xi_{2}$ eigenvector fields, are no longer defined at singularities. Still,
the behavior of tensorlines near a singularity has some analogies,
as well as notable differences, with the behavior of trajectories
of a two-dimensional vector field near fixed point. In the absence
of symmetries, there are two structurally stable singularities of
a tensorline field: trisectors and wedges.

\emph{Trisector} singularities are similar to saddle points in two-dimensional
flows, except that they have three (as opposed to two) distinguished
strainlines asymptotic to them (Fig. \ref{fig:tri_wedge}).

\emph{Wedge} singularities are a mix between a saddle and a source
or a sink. On the one hand, there is a continuous family of infinitely
many neighboring tensorlines asymptotic to a wedge. At the same time,
a wedge also has discrete tensorlines asymptotic to it, resembling
the stable and unstable manifolds of a saddle (Fig. \ref{fig:tri_wedge}).

\subsection{Numerical detection of singularities}

At a singularity in an incompressible flow, the elements of the Cauchy--Green
strain tensor satisfy 
\begin{equation}
C_{11}-C_{22}=0\ \ \mbox{and}\ \ C_{12}=0,\label{eqn:Sing_conditions-1}
\end{equation}
where $C_{ij}$ is the $(i$,$j)$ the entry of $C$. The singularities
are, therefore, precisely points where the zero level-curves of the
scalar functions $f=C_{11}-C_{22}$ and $g=C_{12}$ intersect. These
intersections can be found by linearly interpolating $f$ and $g$
along the edges of a numerical grid \cite{Ten_lines_2D}.

In regions of high mixing and chaos, the entries of the Cauchy--Green
strain tensor can be large and noisy. An indication of noise in an
incompressible flow is that the determinant of $C$ is far from being
equal to $1$. These noisy points result in spurious intersection
of the zero levels of $f$ and $g$, and hence spurious singularity
detection.

A crude but effective way of filtering out most if the these spurious
intersections is to consider only parts of the zero level set of $f$
and $g$ on which $|\lambda_{1}\lambda_{2}-1|>1$ holds.

\subsection{Numerical classification of singularities}
Once the singularities are located, we need a robust procedure to
classify each of these singularities as a wedge or a trisector. The existing methods for distinguishing trisector singularities of a tensor field from its wedge singularities require further differentiation of the tensor field \cite{Tricoche}. In our experience, this introduces further
noise affecting the robustness of the results. Here, we introduce a differentiation-free method for identifying trisectors
and wedges. This method also is used to find the direction of the separatrices emanating from a trisector.

A distinguishing geometric feature
of a trisector singularity is the three separatrices emanating from
it. Close enough to the singularity, these separatrices are close
to straight lines. Therefore, the separatrices will be approximately
perpendicular to a small circle centered at the singularity. Consequently.
the intersection of the trisectors with the circle approximately maximizes
the quantity 

\begin{equation}
f_{i}(\theta)=\frac{|\left<\xi_{i},r\right>|}{\left|\xi_{i}\right|\left|r\right|}\label{eq:rDOTxi}
\end{equation}
associated with the vector field $\xi_{i}$, with $r$ is the vector
from the singularity pointing towards the point $\theta$ on the small
circle.

\begin{figure}[H]
\begin{center}
\includegraphics[width=0.3\textwidth]{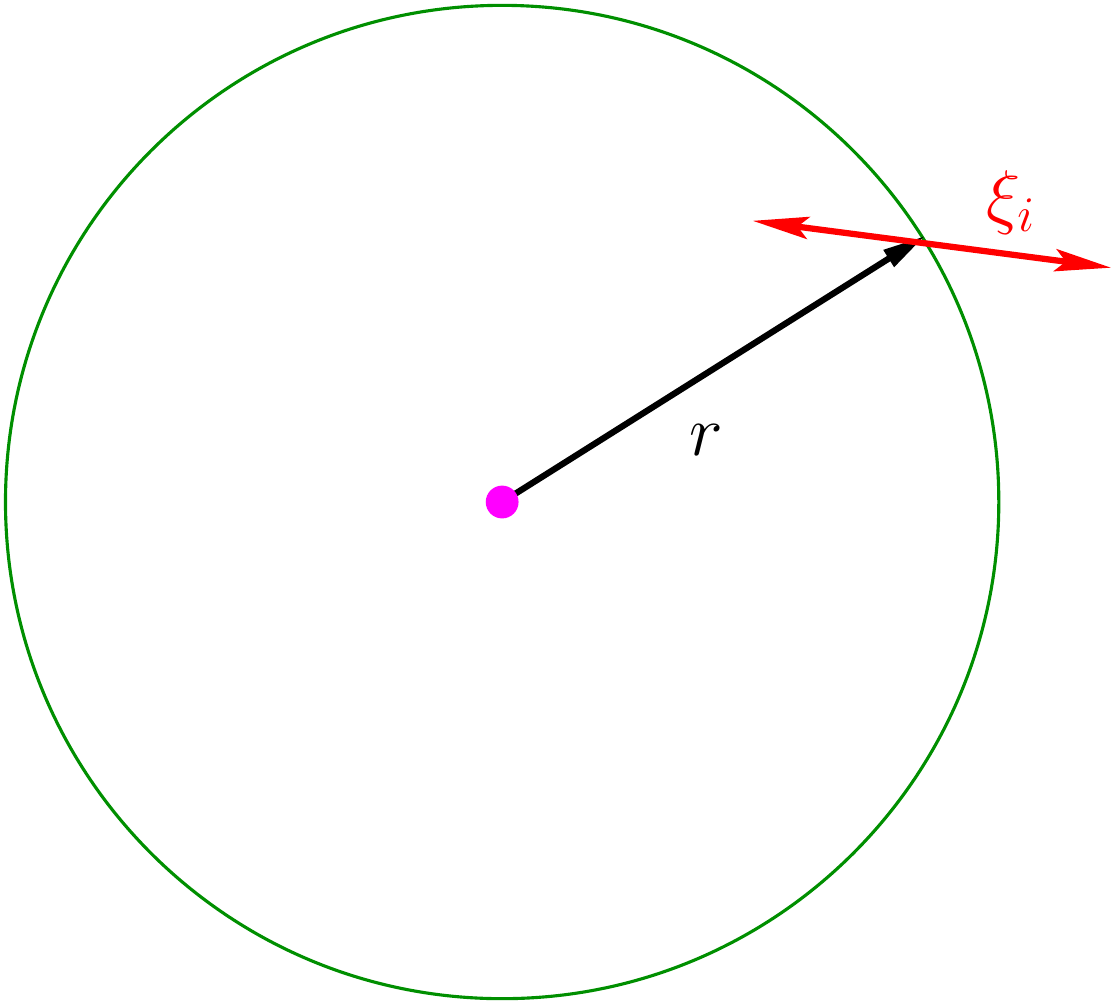} 
\end{center}
\caption{In equation (\ref{eq:rDOTxi}), $f(\theta)$ is defined as the normalized
inner product of $r$ and $\xi_{i}$.}
\label{fig:rDOTxi} 
\end{figure}

For a trisector, $f_i(\theta)$ assumes the value $0$ and $1$ three
times, with $0$'s and $1$'s alternating, as $\theta$ increases
from $0$ to $2\pi$. In contrast, for a wedge, $f_i$ assumes $1$
three times, and assumes a zero value only once. We use this difference
between wedges and trisectors in identifying them numerically.

Moreover, for a trisector, the $\theta$ values for which $f_i(\theta)=1$
indicate the direction of its separatrices corresponding to the vector field $\xi_i$.

\subsection{Structurally stable heteroclinic tensorlines and their numerical
detection}

As a consequence of trisector and wedge geometries, there can be no
unique connection between two wedge singularities. Indeed, if there
is one such connection, there must be infinitely many. On the other
hand, as in the case of heteroclinic orbits between saddles of an
ODE, trisector-trisector connections are structurally unstable. Therefore,
the only types of tensorlines connecting two singularities of the
Cauchy--Green strain tensor that are locally unique and structurally
stable are trisector-wedge connections. 

The numerical detection of trisector-wedge connections proceeds by
tracking the separatrices leaving a trisector, and monitoring whether
they enter the attracting sector of a small circle surrounding a wedge.

\end{appendices}

\end{document}